\documentclass[journal,10pt]{IEEEtran}
\usepackage[colorlinks,linkcolor=black,anchorcolor=black,linktocpage=true,urlcolor=blue,citecolor=black]{hyperref}

\usepackage{stfloats}
\usepackage{booktabs}
\usepackage{algorithmic}
\usepackage{tikz}
\usetikzlibrary{arrows}
\usepackage{array}
\usepackage{graphicx}

\usepackage{url}
\usepackage{cite}
\usepackage{amsmath,amsthm,amsfonts,amssymb,mathrsfs,bm}
\allowdisplaybreaks[4]
\usepackage{tabularx,array,enumerate,pifont,tikz,epstopdf}
\newtheorem{theorem}{Theorem}
\newtheorem{lemma}{Lemma}

\newtheorem{remark}{Remark}

\newtheorem{assumption}{Assumption}

\begin{document}

\title{\hspace*{-2.8mm} \LARGE \bf Optimal Output Consensus of Heterogeneous Linear Multi-Agent Systems Over Weight-Unbalanced Directed Networks}

\author{Jin~Zhang,
        Lu~Liu,~\IEEEmembership{Senior Member, IEEE,}
        Haibo~Ji,
        and~Xinghu~Wang 
\thanks{ Manuscript submitted to \textit{IEEE Transactions on Cybernetics}; May 22, 2020; revised April 23, 2021.
	
	J. Zhang is with the Department of Automation, University of Science and Technology of China, Hefei 230027, China, and also with the Department of Biomedical Engineering, City University of Hong Kong, Hong Kong (e-mail: zj55555@mail.ustc.edu.cn).
	
	L. Liu is with the Department of Biomedical Engineering, City University of Hong Kong, Hong Kong (e-mail: luliu45@cityu.edu.hk).
	
	H. Ji and X. Wang are with the Department of Automation, University of Science and Technology of China, Hefei 230027, China (e-mail: jihb@ustc.edu.cn; xinghuw@ustc.edu.cn).}}


\maketitle

\begin{abstract}
		This paper investigates the distributed optimal output consensus problem of heterogeneous linear multi-agent systems over weight-unbalanced directed networks. A novel distributed continuous-time state feedback controller is proposed to steer the outputs of all the agents to converge to the optimal solution of the global cost function. Under the standard condition that the unbalanced digraph is strongly connected and the local cost functions are strongly convex with global Lipschitz gradients, the exponential convergence of the closed-loop multi-agent system is established. Then, the proposed state feedback control law is extended to an observer-based output feedback setting. Two examples are finally provided to illustrate the effectiveness of the proposed control schemes. 
\end{abstract}


\IEEEpeerreviewmaketitle

\section{Introduction}

\IEEEPARstart{O}{ver} the past decade, the distributed optimization problem (DOP) has attracted increasing attention due to its broad applications in sensor networks, distributed parameter estimation, power systems and machine learning; see, for example, \cite{molzahn2017survey,ram2010distributed}. In the DOP, each agent is assigned with an individual local cost function. The objective is to minimize the sum of all local cost functions in a distributed manner by using only neighboring information and local computation. Seminal works on this topic can be traced back to \cite{tsitsiklis1984problems,nedic2009distributed}, and for recent progress, one may refer to the relevant reviews \cite{yang2019survey,nedic2018distributed} and references therein. 

Most of the existing works focus on discrete-time dynamics, to name a few, see \cite{liang2019distributed,nedic2014distributed,xi2017dextra,xi2016distributed,xi2018linear,wang2019distributed}. However, since many practical systems operate in the continuous-time setting, such as unmanned vehicles and robots, a few efforts have been made recently for the case of continuous-time dynamics \cite{wang2010control,gharesifard2013distributed,kia2015distributed}. Based on the first-order gradients, the authors in \cite{wang2010control} propose a distributed proportional-integral (PI) control scheme for multi-agent systems over undirected graphs. Then the approach is extended to deal with weighted-balanced directed networks in \cite{gharesifard2013distributed}. To eliminate the requirement for additional information communication by PI feedback, the modified Lagrangian based (MLB) algorithm is then developed in \cite{kia2015distributed} at the expense of special initialization.

To handle the DOP on weight-unbalanced directed graphs, some consensus based protocols are commonly used, such as the push-pull based protocols \cite{xi2016distributed,xi2017dextra} and the push-sum based protocols \cite{touri2016saddle,nedic2014distributed}. However, the above-mentioned protocols usually involve certain global information, including in-degree \cite{touri2016saddle} and out-degree \cite{xi2016distributed,nedic2014distributed}, which might not be available in general directed networks. A distributed continuous-time control strategy is designed in \cite{li2017distributed} to deal with the weight-unbalanced directed graphs, but it cannot achieve the optimal solution when the left eigenvector is not available in advance. To handle the imbalance, a distributed discrete-time algorithm is proposed in \cite{xi2018linear} with the gradient being divided by an additional variable, which is designed to exponentially converges to the left eigenvector corresponding to the eigenvalue one of the row-stochastic matrix. More recently, the discrete-time algorithm in \cite{xi2018linear} is extended to a continuous-time version in \cite{zhu2018continuous}, and the explicit dependency on the left eigenvector is removed in comparison with the algorithm in \cite{li2017distributed}. 

It is worth mentioning that the conventional DOP in the aforementioned works can be regarded as a distributed optimal output consensus (OOC) problem with single integrator dynamics. However, there are quite a few engineering tasks in practice that could be reformulated as the OOC problem for more general agent dynamics, such as the economic dispatch in power systems \cite{stegink2016Unifying}, rigid body attitude formation control \cite{song2017Relative} and source seeking in multi-robot systems \cite{zhang2011extremum}. Recently, many efforts are dedicated to solving the OOC for double integrators \cite{tran2018distributed}, high-order linear systems \cite{xie2019global,tang2018optimal} and nonlinear systems \cite{wang2015distributed,tang2018distributed}, over undirected graphs. It can be noted that the control design in such scenarios is much more challenging. Typically, the control design for the OOC problem can be classified into two types. The first type is a control scheme based on the two-layer structure consisting of an optimal signal generator and a reference-tracking controller \cite{tang2018optimal,tang2018distributed}. However, this type of control design requires the agent dynamics to be minimum-phase and have well-defined vector relative degrees \cite{tang2018optimal}, and may fail in the scenario when some of the optimal signal generators result in an augmented system for which the reference-tracking problem cannot be solved. Different from the first type, the second one concentrates on developing integrated control laws, which avoids the requirement on the optimal signal as a reference \cite{zhao2017distributed,li2019distributed}. The authors in \cite{zhao2017distributed} propose two adaptive control laws to address the OOC problem for homogeneous linear multi-agent systems. However, the controllers can only be applied when the gradients satisfy a specific structure, and, certain global information is needed to verify their applicability. In a recent work \cite{li2019distributed}, the OOC problem of heterogeneous linear multi-agent systems is reformulated as a special output regulation problem, which is then solved by designing a new controller based on the solutions of well-designed linear matrix equations. Unfortunately, the aforementioned controllers can only be applied to undirected graphs.

In our preliminary work \cite{zhang2020exponential}, the OOC problem for homogeneous linear multi-agent systems over weight-unbalanced directed graphs is solved by designing a distributed state feedback controller. The design of the state feedback controller requires a prior knowledge of the left eigenvector corresponding to the eigenvalue zero of the Laplacian matrix $\mathcal{L}$. In this paper, we consider the OOC problem over weight-unbalanced directed graphs for heterogeneous linear multi-agent systems. A novel distributed state feedback controller is firstly proposed by introducing an additional variable to avoid the explicit dependence on the left eigenvector. To tackle heterogeneous linear agent dynamics, we take advantage of the well-designed matrix equations from \cite{li2019distributed}, which serves as a modification of regulator equations in \cite{huang2004nonlinear}. Then the proposed state feedback controller is further extended to an observer-based output feedback one. The main contributions of this paper are summarized as follows:

1) This work investigates the OOC problem for general linear multi-agent systems on weight-unbalanced digraphs. Compared with most existing works, the scenarios considered in this paper are much more general and thus more applicable in practice. On one hand, in contrary to integrator-type agent dynamics discussed in \cite{gharesifard2013distributed,zhu2018continuous}, we consider more general heterogeneous linear systems. On the other hand, unlike the works of solving the OOC problem for linear multi-agent systems on undirected networks \cite{zhao2017distributed,tang2018optimal,li2019distributed}, we focus on more general and thus more challenging weight-unbalanced directed networks. To address the challenges arising from the asymmetry of the Laplacian matrices corresponding to directed graphs, we utilize some useful results from Kronecker matrix algebra and the direct sum operation of vectors instead of the commonly used orthogonal transformation.
	
2) When the state information is not measurable, which is often the case in practice, we extend the newly developed state feedback controller to an observer-based output feedback one so that the considered OOC problem can still be solved. Therefore, compared with the existing results in \cite{zhao2017distributed,xie2019global}, our work further expands the scope of distributed optimization in practical applications.

The rest of this paper is organized as follows. Preliminaries and problem formulation are presented in Section \ref{section preliminaries} and \ref{section problem formulation}, respectively. Design of control laws and analysis of the resulting closed-loop systems are provided in Section
\ref{section main results} followed by illustrative examples in Section \ref{section simulation results}. Finally, the conclusion is stated in Section \ref{section conclusion}.

\textit{Notations}: Let $\mathbb{R}^{n}$ and $ \mathbb{R}^{n \times m}$ denote the sets of real vectors of dimension $n$ and real matrices of dimension $n \times m,$ respectively. Let $\mathbf{1}_{n}$ and $\mathbf{0}_{n}$ denote the column vector of dimension $n$ with all entries equal to one and zero, respectively. Let $I_{n}$ denote the identity matrix of dimension $n$. Let $\otimes$ denote the Kronecker product of matrices. For a matrix $A \in \mathbb{R}^{n \times n}$, $A^{\mathrm{T}}$ and $ \operatorname{tr}(A) $ represent its transpose and trace, respectively. $\|\cdot\|$ represents the Euclidean norm of vectors or the induced 2-norm of matrices. $ \operatorname{col}\left(x_{1}, x_{2}, \ldots, x_{n}\right) $ represents a column vector with $x_{1}, x_{2}, \ldots, x_{n}$ being its elements. $\operatorname{diag}\left(x_{1}, x_{2}, \ldots, x_{n}\right)$ represents a diagonal matrix with $x_{1}, x_{2}, \ldots, x_{n}$ being its diagonal elements;  $\operatorname{diag}\left(B_{1},B_{2}, \ldots, B_{N}\right)$ represents a block diagonal matrix with matrices $B_{i} \in \mathbb{R}^{n_{i} \times p_{i}}, i \in 1,2, \ldots N$ being its diagonal block elements. For a differentiable function $f: \mathbb{R}^{n} \rightarrow \mathbb{R}$, $\nabla f$ is its gradient.

\section{Preliminaries}\label{section preliminaries}
In this section, we present some preliminaries on graph theory, convex analysis, Kronecker matrix algebra, and perturbed system theory.

\subsection{Graph Theory} 
A directed graph (in short, a digraph) is adopted to depict the agent information flow. A weighted directed graph of order $ N $ is a triplet $\mathcal{G}=(\mathcal{V}, \mathcal{E}, \mathcal{A})$, where $\mathcal{V}=\{1,2, \ldots, N\}$ is a set with $ N $  vertices called nodes, $\mathcal{E} \subseteq \mathcal{V} \times \mathcal{V}$ is a set of ordered pairs of nodes called edges, and $\mathcal{A}=\left[a_{i j}\right] \in \mathbb{R}^{N \times N}$ is an associated weighted adjacency matrix. For $ i,j\in\mathcal{V} $, the ordered pair $(j, i) \in \mathcal{E}$ denotes an edge from $ j $ to $ i $, that is, the $i$th agent can receive information from the $j$th agent, but not vice versa. In this case, $ j $ is called an in-neighbor of $ i $, and $ i $ is called an out-neighbor of $ j $. A directed path is an ordered sequence of nodes such that any pair of consecutive nodes in the sequence is a directed edge. A digraph is strongly connected if there exists a directed path in each direction between each pair of nodes. The associated adjacency matrix $\mathcal{A}$ is defined as $a_{i j}>0$ if $(j, i) \in \mathcal{E}$, otherwise $a_{i j}=0 $, and $a_{i i}=0$ for all $ i\in \mathcal{V}$ since it is assumed that there are no self-loops in a digraph. Furthermore, the Laplacian matrix $\mathcal{L}=\left[l_{i j}\right] \in \mathbb{R}^{N \times N}$ associated with the digraph $\mathcal{G}$ is defined as $l_{i i}=\sum_{j=1}^{N} a_{i j}$ and $l_{i j}=-a_{i j}$ for $i \neq j $. A digraph $\mathcal{G}$ is called weight balanced iff $\mathbf{1}_{N}^{\mathrm{T}} \mathcal{L}=\mathbf{0}_{N}^{\mathrm{T}}$. For a more detailed introduction of graph theory, please refer to \cite{bullo2019lectures}.

\begin{lemma}\cite{olfati-Saber2004consensus,ren2005consensus}
	Assume that the unbalanced directed graph $ \mathcal{G} $ is strongly connected. Let $ \mathcal{L} $ be the associated Laplacian matrix. Then
	\begin{itemize}
		\item there exists a positive left eigenvector $ r=\left(r_{1},r_{2},\ldots,r_{N} \right)^{\mathrm{T}}  $ associated with the eigenvalue zero such that $ r^{\mathrm{T}}\mathcal{L}=\mathbf{0}_{N}^{\mathrm{T}} $ and $ \sum_{i=1}^{N}r_{i}=1 $;
		\item $\bar{\mathcal{L}}=\left( R\mathcal{L}+\mathcal{L}^{\mathrm{T}}R\right) /2 $ is positive semidefinite, where $ R=\mathrm{diag}\left(r_{1},r_{2},\ldots,r_{N} \right) $, and the eigenvalues of $ \bar{\mathcal{L}} $ can be ordered as $ 0=\lambda_{1}<\lambda_{2}\leq \lambda_{3}\leq\ldots \leq\lambda_{N} $; and
		\item $ \exp\left( -\mathcal{L}t \right)  $ is a nonnegative matrix with positive
		diagonal entries for all $ t>0 $, and $ \lim_{t\to\infty}\exp\left( -\mathcal{L}t \right)= \mathbf{1}_{N}r^{\mathrm{T}} $.
		\label{graph theory lemma}
	\end{itemize}
	
\end{lemma}

\subsection{Convex Analysis}
In this subsection, the definitions of strong convexity and Lipschitz continuity are given, see  \cite{bertsekas2009convex} for more details.

A continuously differentiable function $f: \mathbb{R}^{n} \rightarrow \mathbb{R}$ is strongly convex on $\mathbb{R}^{n}$ if there exists a positive constant $m $ such that $(x-y)^{T}(\nabla f(x)-\nabla f(y)) \geq m\|x-y\|^{2}$ for all $x, y \in \mathbb{R}^{n} $. A function $g: \mathbb{R}^{n} \rightarrow \mathbb{R}^{n}$ is globally Lipschitz on $\mathbb{R}^{n}$ if there exists a positive constant $M $ such that $\|g(x)-g(y)\| \leq M\|x-y\|$ for all $x, y \in \mathbb{R}^{n}$.

\subsection{Kronecker Matrix Algebra}
For $ E\in \mathbb{R}^{n\times m} $, Let $ \mathrm{col}_{i}(E) \in \mathbb{R}^{n} $ be the $ i $th column of matrix $ E $, and denote
\begin{equation}
\operatorname{vec}(E) \triangleq\left(\begin{array}{c}{\operatorname{col}_{1}(E)} \\ {\vdots} \\ {\operatorname{col}_{m}(E)}\end{array}\right) \in \mathbb{R}^{n m}\notag
\end{equation}
as the column vector of dimension $ nm $ obtained by stacking $ \mathrm{col}_{i}(E) $. The following two crucial lemmas are given in \cite{bernstein2009matrix}.

\begin{lemma}\label{lemma 2}
	Let $ E\in \mathbb{R}^{n\times m} $ and $ F\in \mathbb{R}^{m\times n} $. Then, $\operatorname{tr}(E F)= \operatorname{tr}(F E)=\left(\operatorname{vec}\left(E^{\mathrm{T}}\right)\right)^{\mathrm{T}} \operatorname{vec}(F) $. Additionally, for real column vectors $ a\in \mathbb{R}^{n} $ and $ b\in \mathbb{R}^{n} $, $\operatorname{tr}\left(b a^{T}\right)= a^{T} b$.
	
\end{lemma}

\begin{lemma} \label{lemma 3}
	Let $ E\in \mathbb{R}^{n\times m} , F\in \mathbb{R}^{m\times l} $, and $ G\in \mathbb{R}^{l\times k} $. Then, $\operatorname{vec}(E F G)=\left(G^{\mathrm{T}} \otimes E\right) \operatorname{vec}(F)$.
\end{lemma}

\subsection{Perturbed System Theory}
Last but not least, the theory of perturbed systems which facilitates subsequent analysis is addressed.

\begin{lemma} \label{lemma 4}
	Consider the perturbed system
	\begin{equation} \label{lemma_perturbed system}
		\dot{x}=g(t,x)+\kappa(t,x).
	\end{equation}
	Let $x=0$ be an exponentially stable equilibrium point of the nominal system $ \dot{x}=g(t,x) $, where $ g $ is continuously differentiable and the Jacobian matrix $ [\partial g/ \partial x] $ is bounded on $ \mathbb{R}^{n} $. Suppose the perturbation term $\kappa(t, x)$ satisfies $ \kappa(t,0)=0 $ and $ \|\kappa(t, x)\| \leq \gamma(t)\|x\| $, where $\lim_{t\to\infty}\gamma(t)= 0$. Then, the origin is an exponentially stable equilibrium point of the perturbed system (\ref{lemma_perturbed system}). 
	
	\begin{proof}
		It can be proved by a simple combination of Corollary 9.1 and Lemma 9.5 in \cite{khalil2002nonlinear}.
	\end{proof}
\end{lemma}

\section{Problem Formulation} \label{section problem formulation}

Consider a multi-agent system of $ N $ heterogeneous agents described by the following dynamics,
\begin{equation}
\begin{aligned}
\dot{x}_{i} & = A_{i} x_{i}+B_{i} u_{i}, \\
y_{i} & = C_{i} x_{i}, \quad i=1,2, \ldots, N,
\end{aligned}\label{dynamics}
\end{equation}
where $ x_{i}\in \mathbb{R}^{n_{i}} $, $ u_{i}\in \mathbb{R}^{p_{i}} $ and $ y_{i}\in \mathbb{R}^{q} $ are the state, control input and output of the $ i $th agent, respectively. $ A_{i}\in \mathbb{R}^{n_{i} \times n_{i}} $, $ B_{i}\in \mathbb{R}^{n_{i} \times p_{i}} $, and $ C_{i}\in \mathbb{R}^{q \times n_{i}} $ are constant matrices.

A local cost function $ f_{i}(y): \mathbb{R}^{q}\to \mathbb{R} $ is assigned to each agent $ i $ for $ i=1,2,\ldots,N $, which is only available to agent $ i $. Define the global cost function as $ f(y)=\sum_{i=1}^{N}f_{i}(y) $. The objective of this work is to design distributed controllers such that the outputs of all the agents are steered to the optimal solution $ y^{*} $ of the following optimization problem,
\begin{equation}
\min_{y\in \mathbb{R}^{q}} f(y). \label{objetive}
\end{equation}

\begin{remark}
	Unlike the output consensus problem of heterogeneous linear multi-agent systems \cite{wieland2011internal,kim2010output}, we consider the more general and difficult optimal output consensus problem in this paper. The differences and also challenges in this scenario are to steer the outputs of all the agents not only to achieve consensus but also to reach the optimal solution of the global cost function. In this case, the design of distributed controllers is more challenging.
\end{remark}

To achieve the objective, we need some standard assumptions.

\begin{assumption}\label{cost function assumption}
	For $ i=1,2,\ldots,N $, the local cost function $ f_{i} $ is continuously differentiable and strongly convex with constant $ m_{i} $, and $\nabla f_{i}$ is globally Lipschitz on $ \mathbb{R}^{q} $ with constant $ M_{i} $.
	
\end{assumption}

\begin{remark}
	The strong convexities of the local cost functions in Assumption \ref{cost function assumption} guarantee that the optimal solution $ y^{*}\in \mathbb{R}^{q} $ is existing and unique. Assumption \ref{cost function assumption} is standard and thus commonly used in many existing works, see \cite{kia2015distributed,zhu2018continuous}.
\end{remark}

\begin{assumption} 	\label{graph assumption}
	The communication directed graph $ \mathcal{G} $ is strongly connected.

\end{assumption}

Define $ \tilde{f}(Y)=\sum_{i=1}^{N}f_{i}(y_{i}) $ with $ Y=\mathrm{col}(y_{1}, y_{2},\ldots,y_{N}) \in \mathbb{R}^{qN} $. Then similar to the previous works \cite{li2017distributed,li2019distributed}, under Assumption \ref{graph assumption}, we can reformulate problem (\ref{objetive}) as
\begin{equation}
\begin{aligned}
&\min _{Y \in \mathbb{R}^{qN}}  \tilde{f}(Y), \quad  \tilde{f}(Y) = \sum_{i=1}^{N} f_{i}\left(y_{i}\right), \\
&\text { subject to }  \quad \left(\mathcal{L} \otimes I_{q}\right) Y =\mathbf{0}.
\end{aligned} \label{equivalent problem}
\end{equation}

To handle linear agent dynamics, we need an additional assumption.
\begin{assumption}
	For $ i=1,2, \ldots, N $, $ \left( A_{i}, B_{i} \right) $ is stabilizable and
	\begin{equation}
	\operatorname{rank}\left[\begin{array}{cc}
	{C_{i} B_{i}} & {\mathbf{0}_{q \times p_{i}}} \\
	{-A_{i} B_{i}} & {B_{i}}
	\end{array}\right]=n_{i}+q. 
	\notag 
	\end{equation} 
	\label{matrix equations assumption}
\end{assumption}

\begin{lemma} \cite{li2019distributed} \label{lemma_matrix equations}
	Suppose that Assumption \ref{matrix equations assumption} holds. Then for $ i=1,2, \ldots, N $, the linear matrix equations
	\begin{subequations}	
		\begin{align}	
		C_{i} \Psi_{i}-I_{q} &=\mathbf{0}_{q \times q}, \label{matrix equations_Psi} \\ 
		B_{i} \Phi_{i}-A_{i} \Psi_{i} &=\mathbf{0}_{n_{i} \times q}, \label{matrix equations_Phi} \\ 
		B_{i} \Upsilon_{i}-\Psi_{i} &=\mathbf{0}_{n_{i} \times q}  \label{matrix equations_Upsilon}
		\end{align}	\label{matrix equations}
	\end{subequations} 
have solution triplets $ \left( \Upsilon_{i}, \Phi_{i}, \Psi_{i} \right) $.
\end{lemma}

\begin{remark}
	Although it is proposed in \cite{tang2018optimal} to handle the OOC problem for linear systems having well-defined vector relative degrees, the two-layer control scheme would fail to achieve the objective if the agent dynamics are non-minimum phase. This would thus limit the scope of the algorithm's application as many practical systems are non-minimum phase. On the contrary, our approach in this work does not suffer such a restriction.  
\end{remark}

\section{Main Results}\label{section main results}

In this section, a distributed state feedback controller is firstly developed to deal with the OOC problem over weight-unbalanced digraphs. Then the proposed distributed state feedback controller is extended to a distributed observer-based output feedback controller when the state information is not measurable.

\subsection{Distributed state feedback controller}

In this subsection, we design the distributed state feedback controller. Based on the relative outputs and individual agent state, the controller for each agent is designed as follows,
\begin{subequations}\label{control law}	
	\begin{align}	
	u_{i} & =-K_{i}x_{i}+\Upsilon_{i}\omega_{i}-\left(\Phi_{i}-K_{i}\Psi_{i} \right) \rho_{i},\label{control law a} \\ 
	\dot{\rho}_{i} & =\omega_{i} :=-\dfrac{\nabla f_{i}(y_{i})}{z_{i}^{i}} -\gamma_{1}\sum_{j=1}^{N}a_{i j}\left(y_{i}-y_{j} \right)-\gamma_{2} v_{i}, \label{control law b} \\ 
	\dot{v}_{i} & =\gamma_{1} \sum_{j=1}^{N} a_{i j}\left(y_{i}-y_{j} \right), \label{control law c} \\ 
	\dot{z}_{i} & =-\sum_{j=1}^{N} a_{i j}\left(z_{i}-z_{j} \right) \label{control law d} ,
	\end{align}	
\end{subequations}
where $ z_{i}\in \mathbb{R}^{N},i=1,2,\ldots,N $, with $ z_{i}^{i} $ being its $ i $th component and with its initial value $  z_{i}(0) $ satisfying $ z_{i}^{i}(0)=1 $, $ z_{i}^{j}(0)=0 $ for all $ j\neq i $. $ \rho_{i}\in \mathbb{R}^{q} $ and $ v_{i}\in \mathbb{R}^{q} $ are auxiliary variables with initial value $ v_{i}(0)=0 $, $ \omega_{i}\in \mathbb{R}^{q} $ is an intermediate state. $ \gamma_{1} $ and $ \gamma_{2} $ are positive constants, $ K_{i}\in\mathbb{R}^{p_{i}\times n_{i}} $ is a feedback gain matrix such that $ A_{i}-B_{i}K_{i} $ is Hurwitz, and $ (\Upsilon_{i},\Phi_{i},\Psi_{i}) $ are the solution triplets of linear matrix equations (\ref{matrix equations}). It is worth noting that control law (\ref{control law a}) is composed of the relative outputs and their first-order and second-order integrals. Equations (\ref{control law b}) and (\ref{control law c}) can be seen as a modification to the well-known MLB algorithm in \cite{kia2015distributed} if setting $ \rho_{i}=y_{i} $ and $ z_{i}^{i}=1 $. 

Let $ x=\mathrm{col}(x_{1},x_{2},\ldots,x_{N}) $, $ \rho=\mathrm{col}(\rho_{1},\rho_{2},\ldots,\rho_{N}) $, $ v=\mathrm{col}(v_{1},v_{2},\ldots,v_{N}) $, $ z=\mathrm{col}(z_{1},z_{2},\ldots,z_{N}) $, $ Z_{N}=\mathrm{diag}(z_{1}^{1},z_{2}^{2},\ldots,z_{N}^{N}) $, $ A=\mathrm{diag}(A_{1},A_{2},\ldots,A_{N}) $, $ B=\mathrm{diag}(B_{1},B_{2},\ldots,B_{N}) $, $ C=\mathrm{diag}(C_{1},C_{2},\ldots,C_{N}) $, $ K=\mathrm{diag}(K_{1},K_{2},\ldots,K_{N}) $, $ \Upsilon=\mathrm{diag}(\Upsilon_{1},\Upsilon_{2},\ldots,\Upsilon_{N}) $, $ \Phi=\mathrm{diag}(\Phi_{1},\Phi_{2},\ldots,\Phi_{N}) $, $ \Psi=\mathrm{diag}(\Psi_{1},\Psi_{2},\ldots,\Psi_{N}) $, and $ \nabla \tilde{f}(Y)= \mathrm{col}\big( \nabla f_{1}(y_{1}),\nabla f_{2}(y_{2}),\ldots, \nabla f_{N}(y_{N}) \big)  $. Then, by substituting the control law (\ref{control law}) into the dynamics (\ref{dynamics}), the closed-loop system can be written in the following compact form,
\begin{subequations}	\label{closed-loop system}
	\begin{align}	
	\dot{x} & =\left( A-BK \right)x+B\Upsilon \dot{\rho}-\left( B\Phi-BK\Psi \right) \rho ,\label{closed-loop system_a} \\ 
	\dot{\rho} & =-\left( Z_{N}^{-1}\otimes I_{q} \right)\nabla\tilde{f}(Y)-\gamma_{1} \left( \mathcal{L}\otimes I_{q}\right)Y-\gamma_{2} v  , \label{closed-loop system_b}\\
	\dot{v} & =\gamma_{1} \left( \mathcal{L}\otimes I_{q}\right)Y,  \label{closed-loop system_c} \\
	\dot{z} & =-\left(\mathcal{L}\otimes I_{N} \right)z. \label{closed-loop system_d} 
	\end{align}	
\end{subequations}

\begin{remark} 
	Under Asumption \ref{graph assumption}, one can obtain from Lemma \ref{graph theory lemma} that $ z_{i}^{i}(t)>0 $ for all $ t \geq 0 $, which indicates that $ Z_{N}^{-1} $ is well defined \cite{zhu2018continuous}. The term $ \left( Z_{N}^{-1}\otimes I_{q} \right)\nabla\tilde{f}(Y) $ in (\ref{closed-loop system_b}) is used to tackle the imbalance caused by employing only asymmetric Laplacian matrix $ \mathcal{L} $, where $ Z_{N}^{-1} $ serves as the role of $ R^{-1} $ in the algorithm in \cite{li2017distributed}. Thus, the purpose of introducing (\ref{closed-loop system_d}) is to estimate the left eigenvector $ r $ associated with the eigenvalue zero of $ \mathcal{L} $ to reduce the restrictive requirement on global information.
	
\end{remark}

To proceed, we first show that the matrix $ Z_{N}^{-1} $ will exponentially tend to $ R^{-1} $. By using Lemma \ref{graph theory lemma}, one can obtain that $ \lim_{t\to\infty}z(t)=\lim_{t\to\infty}\exp\left( -(\mathcal{L}\otimes I_{N})t\right)z(0)=  \left( \mathbf{1}_{N}r^{\mathrm{T}}\otimes I_{N}\right)z(0)=\mathbf{1}_{N}\otimes r  $. This implies that $ \lim_{t\to\infty} Z_{N}^{-1} = R^{-1} $ exponentially. 

Hereinafter, to cope with the difficulties generated by asymmetric information transmission, we utilize the theories of perturbed systems and input-to-state stability. Define a new variable $ \xi=\mathrm{col}\left(x,\rho,v \right)  $. Then, (\ref{closed-loop system_a})-(\ref{closed-loop system_c}) can be rewritten as follows,
\begin{align}\label{perturbed form}
\underbrace{
	\left(\begin{array}{c}
	{\dot{x}} \\
	{\dot{\rho}}\\
	{\dot{v}}
	\end{array}\right)
}_{\dot{\xi}}
&=\underbrace{
	\left(\begin{array}{c}
	{\left( A-BK \right)x+B\Upsilon \dot{\rho}-\left( B\Phi-BK\Psi \right) \rho}\\
	{-\left( R^{-1}\otimes I_{q} \right)\nabla\tilde{f}(Y)-\gamma_{1} \left( \mathcal{L}\otimes I_{q}\right)Y-\gamma_{2} v} \\
	{\gamma_{1} \left( \mathcal{L}\otimes I_{q}\right)Y}
	\end{array}\right)
}_{g(\xi)} \notag\\
&+\underbrace{
	\left(\begin{array}{c}
	{\mathbf{0}_{Nq}}\\
	{\left( \left( R^{-1}-Z_{N}^{-1}\right) \otimes I_{q} \right)\big( \nabla\tilde{f}(Y)-\nabla\tilde{f}(\bar{Y})\big) } \\
	{\mathbf{0}_{Nq}}
	\end{array}\right)
}_{\kappa(t,\xi)} \notag\\ 
&+\underbrace{
	\left(\begin{array}{c}
	{\mathbf{0}_{Nq}}\\
	{\left( \left( R^{-1}-Z_{N}^{-1}\right) \otimes I_{q} \right) \nabla\tilde{f}(\bar{Y}) } \\
	{\mathbf{0}_{Nq}}
	\end{array}\right)
}_{\omega(t)},
\end{align}
where $ \bar{Y}=C\bar{x} $, with $ \bar{x} $ being the component of the equilibrium point $ \bar{\xi}=\operatorname{col}(\bar{x},\bar{\rho},\bar{v}) $ of the following system, 
\begin{equation}\label{subsystem}
\dot{\xi}=g(\xi).
\end{equation} 

In what follows, our primary goal is to show that $ \bar{Y}= \mathrm{col}\left( \bar{y}_{1},\bar{y}_{2},\ldots,\bar{y}_{N} \right)  $ is the optimal solution of problem (\ref{equivalent problem}) and the equilibrium point of the closed-loop system (\ref{closed-loop system}) is exponentially stable.

\vspace{1em}
\begin{lemma}\label{Lemma equilibrium point}
	Consider system (\ref{subsystem}) and suppose that Assumptions \ref{cost function assumption}-\ref{matrix equations assumption} hold. Then $ \bar{Y} $ is the optimal solution of problem (\ref{equivalent problem}).	
	\begin{proof}
		In light of (\ref{matrix equations_Phi}) and (\ref{subsystem}), the point $ \bar{\xi}=\operatorname{col}(\bar{x},\bar{\rho},\bar{v}) $ satisfies
		\begin{subequations}\label{equilibrium point equations}	
			\begin{align}	
			\mathbf{0} & =\left( A-BK \right)\left(\bar{x}-\Psi \bar{\rho} \right)  ,\label{equilibrium_x}\\
			\mathbf{0} & =-\left( R^{-1}\otimes I_{q} \right)\nabla\tilde{f}(\bar{Y})-\gamma_{2}\bar{v}  , \label{equilibrium_rho}\\
			\mathbf{0} & =\gamma_{1} \left( \mathcal{L}\otimes I_{q}\right)\bar{Y}.\label{equilibrium_nu}  
			\end{align}	
		\end{subequations}
		
		It follows from (\ref{equilibrium_nu}) that the vector  $ \bar{Y} $ belongs to the null-space of $ \mathcal{L}\otimes I_{q} $. Thus, we have $ \bar{Y}=1_{N}\otimes \tau $ for some constant vector $ \tau \in \mathbb{R}^{q} $. Left multiplying $ \dot{v}=\gamma_{1} \left( \mathcal{L}\otimes I_{q}\right)Y $ by $ r^{\mathrm{T}} \otimes I_{q} $ leads to $ \left( r^{\mathrm{T}} \otimes I_{q}\right) \dot{v}= \mathbf{0} $. Then, with the assumption that $ v(0)= \mathbf{0} $, one obtains $ \left( r^{\mathrm{T}} \otimes I_{q}\right) \bar{v}= \mathbf{0} $ . Next, left multiplying (\ref{equilibrium_rho}) by $ r^{\mathrm{T}}\otimes I_{q} $ results in $ \left( \mathbf{1}_{N}^{\mathrm{T}}\otimes I_{q}\right) \nabla\tilde{f}(\bar{Y})= \mathbf{0} $, which is equivalent to
		\begin{equation}\label{optimal condition}
		\sum_{i=1}^{N}\nabla f_{i}(\bar{y}_{i})=\mathbf{0}.
		\end{equation}
		
		Replacing $ \bar{y}_{i} $ by $ \tau $ in (\ref{optimal condition}) results in $ \sum_{i=1}^{N}\nabla f_{i}(\tau)=\mathbf{0} $. Note that the optimality condition is $ \sum_{i=1}^{N}\nabla f_{i}(y^{*})=\mathbf{0} $. One thus has $ \tau=y^{*} $. Therefore, it can be concluded that $ \bar{Y}=1_{N}\otimes y^{*} $ is the optimal solution of problem (\ref{equivalent problem}).
	\end{proof}
\end{lemma}

\vspace{1em}
\begin{remark}
	In essence, we only need to meet the requirement of $  \left( r^{\mathrm{T}} \otimes I_{q}\right) v(0)= \mathbf{0}  $ for the initial value of $ v(0) $. However, since $ v $ is an internal variable, we can directly set $ v(0)= \mathbf{0} $ for simplicity.
\end{remark}
\newpage

The first main result of this paper is presented below.
\begin{theorem}\label{Theorem 1}
	Consider system (\ref{dynamics}) and suppose Assumptions \ref{cost function assumption}-\ref{matrix equations assumption} hold. Then there exist positive constants  $ \gamma_{1}$ and $\gamma_{2} $ such that the output $ Y $ exponentially converges to the optimal value $ \bar{Y}= \mathbf{1}_{N}\otimes y^{*} $ of problem (\ref{equivalent problem}) under the state feedback controller (\ref{control law}), with $ y^{*} $ being the optimal solution to problem (\ref{objetive}).
	
	\begin{remark}
		The main feature of this paper is that both general high-order dynamics and weight-unbalanced directed graphs can be handled via our proposed controllers. On one hand, different from existing works that require the agent dynamics to be integrator-type \cite{gharesifard2013distributed,zhu2018continuous} or high-order but minimum phase systems \cite{tran2018distributed,xie2019global,tang2018optimal}, the general dynamics considered in this work are allowed to be non-minimum phase. On the other hand, in contrary to most existing works where undirected graphs \cite{liang2019distributed,wang2010control} or balanced directed graphs \cite{gharesifard2013distributed,kia2015distributed} are considered, the proposed controller (\ref{control law}) is able to tackle weight-unbalanced directed graphs without a prior knowledge of the left eigenvector corresponding to the eigenvalue zero of the Laplacian matrix $\mathcal{L}$. It is worth noting that the left eigenvector is a kind of global information, which is required in \cite{li2017distributed}, but may be unavailable in practical applications.
	\end{remark}
	
	\begin{remark}
		Compared with asymptotic convergence, exponential convergence has many advantages from the perspective of stability analysis and control synthesis. Different from existing works \cite{li2019distributed,tang2018optimal,gharesifard2013distributed,zhu2018continuous} where only asymptotic convergence is obtained, the much more desirable exponential convergence can be established in this paper.
	\end{remark}

	\begin{proof}
	The proof can be accomplished by the following three steps. 
	
	Step 1: The exponential stability of system (\ref{subsystem}) is established.
	Introduce the variables substitution $ \tilde{x}=x-\bar{x} $, $ \tilde{\rho}=\rho-\bar{\rho} $ and $ \tilde{v}=v-\bar{v} $ such that the new equilibrium point is transferred to the origin. Note that (\ref{matrix equations_Upsilon}) and $ Y=Cx $ hold. Thus, the dynamics of $ \tilde{x},\tilde{\rho} $ and $ \tilde{v} $ satisfy
	\begin{subequations}\label{change system}	
		\begin{align}	
		\dot{\tilde{x}} & =A_{c} \left(\tilde{x}-\Psi \tilde{\rho} \right) +\Psi \dot{\tilde{\rho}} ,  \label{dot_tilde_x} \\ 
		\dot{\tilde{\rho}} & =-\left( R^{-1}\otimes I_{q} \right)h-\gamma_{1} \left(\mathcal{L}\otimes I_{q} \right)C\tilde{x}  -\gamma_{2}\tilde{v}  , \\
		\dot{\tilde{v}} & =\gamma_{1} \left( \mathcal{L}\otimes I_{q}\right)C\tilde{x},  
		\end{align}	
	\end{subequations}
	where $ h=\nabla\tilde{f}\left( C(\bar{x}+\tilde{x}) \right)- \nabla\tilde{f}\left( C\bar{x} \right) $, and $ A_{c}=A-BK $ is Hurwitz.
		
	Consider the following positive definite function,
	\begin{equation*}
	\begin{aligned}
		V_{1} = & \dfrac{1}{2}\left( C\tilde{x} \right)^{\mathrm{T}} \left( R\otimes I_{q} \right) \left( C\tilde{x} \right) \\
		&+ \dfrac{1}{2}\left( C\tilde{x} + \tilde{v} \right)^{\mathrm{T}} \left( R\otimes I_{q} \right) \left( C\tilde{x} + \tilde{v} \right).
		\end{aligned}
	\end{equation*}
	By using (\ref{matrix equations_Psi}), the derivative of $ V_{1} $ along (\ref{change system}) is given by
	\begin{align} \label{dot_V_1_1}
	\dot{V_{1}}= & -2\left( C\tilde{x} \right)^{\mathrm{T}}h- \gamma_{1} \left( C\tilde{x} \right)^{\mathrm{T}} \left( R\mathcal{L}\otimes I_{q} \right) \left( C\tilde{x} \right)  \notag\\
	& + 2 \left( C\tilde{x} \right)^{\mathrm{T}} \left( R\otimes I_{q} \right) CA_{c}\left( \tilde{x}-\Psi\tilde{\rho} \right) -\tilde{v}^{^{\mathrm{T}}}h \notag\\
	& -2\gamma_{2}\left( C\tilde{x} \right)^{\mathrm{T}} \left( R\otimes I_{q} \right) \tilde{v} -\gamma_{2} \tilde{v}^{\mathrm{T}} \left( R\otimes I_{q} \right) \tilde{v} \notag\\
	& + \tilde{v}^{\mathrm{T}} \left( R\otimes I_{q} \right) CA_{c}\left( \tilde{x}-\Psi\tilde{\rho} \right).
	\end{align}
	
	Using Lemma \ref{graph theory lemma}, we have $ \left( C\tilde{x} \right)^{\mathrm{T}} \left( R\mathcal{L}\otimes I_{q} \right) \left( C\tilde{x} \right) = \left( C\tilde{x} \right)^{\mathrm{T}} \left( \bar{\mathcal{L}}\otimes I_{q} \right) \left( C\tilde{x} \right) $. By referring to the definition of $ R $ in Lemma \ref{graph theory lemma}, one has $ \tilde{v}^{\mathrm{T}} \left( R\otimes I_{q} \right) \tilde{v} \geq r_{\min} \|\tilde{v}\|^{2} $, where $ r_{\min} = \min{\left\lbrace r_{1},r_{2},\ldots,r_{N}\right\rbrace } $. Let $ M=\max{\left\lbrace M_{1},M_{2},\ldots,M_{N}\right\rbrace } $ and $ m=\min{\left\lbrace m_{1},m_{2},\ldots,m_{N}\right\rbrace } $ respectively. One then has $ \|h\| \leq M \| C\tilde{x} \| $ and $ \left( C\tilde{x} \right)^{\mathrm{T}}h \geq m\lVert C\tilde{x} \rVert^{2} $ by Assumption \ref{cost function assumption}. Note that the inequality $ a^{\mathrm{T}}b \leq \frac{1}{\delta}\lVert a\rVert^{2} +\frac{\delta}{4}\lVert b\rVert^{2} $ is always tenable for any $ \delta > 0 $. Thus, the following inequality is satisfied,
	\begin{equation}
	-\tilde{v}^{^{\mathrm{T}}}h  \leq \dfrac{1}{\delta}\|h\|^{2} + \dfrac{\delta}{4}\|\tilde{v}\|^{2} \leq \dfrac{M^{2}}{\delta} \|C\tilde{x}\|^{2} + \dfrac{\delta}{4}\|\tilde{v}\|^{2}.
	\end{equation}
	With the obtained facts, (\ref{dot_V_1_1}) can be rewritten as follows,
	\begin{align}\label{inequality_1}
	\dot{V_{1}}\leq & -\Big( 2m-\dfrac{M^{2}}{\delta} \Big) \lVert C\tilde{x} \rVert^{2}- \Big(\gamma_{2} r_{\min}-\dfrac{\delta}{4} \Big) \|\tilde{v}\|^{2}  \notag \\
	& -\gamma_{1} \left( C\tilde{x} \right)^{\mathrm{T}} \left( \bar{\mathcal{L}}\otimes I_{q} \right) \left( C\tilde{x} \right) -2\gamma_{2}\left( C\tilde{x} \right)^{\mathrm{T}} \left( R\otimes I_{q} \right) \tilde{v} \notag \\
	& + 2 \left( C\tilde{x} \right)^{\mathrm{T}} \left( R\otimes I_{q} \right) CA_{c}\left( \tilde{x}-\Psi\tilde{\rho} \right) \notag \\
	& + \tilde{v}^{\mathrm{T}} \left( R\otimes I_{q} \right) CA_{c}\left( \tilde{x}-\Psi\tilde{\rho} \right).
	\end{align}

	Let $ 0=\lambda_{1}<\lambda_{2}\leq \lambda_{3}\leq\cdots \leq \lambda_{N} $ denote the ordered eigenvalues of the Laplacian matrix $ \bar{\mathcal{L}} $. There exist orthonormal vectors $ \mathbf{1}_{N} $ and $ \eta_{i},i=2,3,\ldots,N $, such that $ \bar{\mathcal{L}}\mathbf{1}_{N}=\mathbf{0} $ and $ \bar{\mathcal{L}}\eta_{i}=\lambda_{i}\eta_{i} $. Define $ \varLambda=\left( C_{1}\tilde{x}_{1},\ldots, C_{N}\tilde{x}_{N} \right) \in \mathbb{R}^{q\times N} $ and construct $ D=\left( d_{1},d_{2},\ldots,d_{N} \right)\in \mathbb{R}^{q\times N} $ with vectors $ d_{1}= \varLambda \mathbf{1}_{N} \in \mathbb{R}^{q}$, $ d_{i}= \varLambda \eta_{i} \in \mathbb{R}^{q}$, $i=2,3,\ldots,N $. Then, it can be verified that $ D=\varLambda Q $, where $ Q=\left(\mathbf{1}_{N},\eta_{2},\ldots,\eta_{N} \right)\in \mathbb{R}^{N\times N} $. Noticing that $ Q $ is an orthogonal matrix, we have $ \varLambda =D Q^{\mathrm{T}} $. By using Lemma \ref{lemma 3}, one then has $ C\tilde{x}=\mathrm{vec}(\varLambda) = \mathrm{vec}\left(DQ^{\mathrm{T}} \right) = \left(I_{N}\otimes D \right)\mathrm{vec}(Q^{\mathrm{T}}) $. Thus, by applying Lemma \ref{lemma 2} and Lemma \ref{lemma 3}, it follows that  
	\begin{align} \label{inequality_1.1}
	&(C\tilde{x})^{\mathrm{T}} \left(\bar{\mathcal{L}} \otimes I_{q}\right) (C\tilde{x}) \notag\\ 
	&=\big(\operatorname{vec}\big(Q^{\mathrm{T}}\big)\big)^{\mathrm{T}}\big(\bar{\mathcal{L}} \otimes D^{\mathrm{T}} D\big) \operatorname{vec}\big(Q^{\mathrm{T}}\big) \notag\\
	&=\big(\operatorname{vec}\big(Q^{\mathrm{T}}\big)\big)^{\mathrm{T}} \operatorname{vec}\big(D^{\mathrm{T}} D Q^{\mathrm{T}} \bar{\mathcal{L}}\big) \notag\\
	&=\operatorname{tr}\big(Q D^{\mathrm{T}} D Q^{\mathrm{T}} \bar{\mathcal{L}}\big)=\operatorname{tr}\big(D Q^{\mathrm{T}} \bar{\mathcal{L}} Q D^{\mathrm{T}}\big) \notag\\
	&=\sum_{i=2}^{N} \eta_{i}^{\mathrm{T}} \bar{\mathcal{L}} \eta_{i} d_{i}^{\mathrm{T}} d_{i} = \sum_{i=2}^{N} \lambda_{i} \eta_{i}^{\mathrm{T}} \eta_{i} d_{i}^{\mathrm{T}} d_{i} \notag\\
	& \geq \sum_{i=2}^{N} \lambda_{2} \eta_{i}^{\mathrm{T}} \eta_{i} d_{i}^{\mathrm{T}} d_{i} 
	=\lambda_{2} \sum_{i=2}^{N} d_{i}^{\mathrm{T}} d_{i} =\lambda_{2}\|s\|^{2},
	\end{align}
	where $ s=\|\operatorname{col}(d_{2},d_{3},\ldots,d_{N})\| $.
	Similarly, define $ \varPi= \left(\tilde{v}_{1}, \tilde{v}_{2},\ldots,\tilde{v}_{N} \right) \in \mathbb{R}^{q\times N} $. One then has
	\begin{align} \label{Cx_v_1}
	&\left( C\tilde{x} \right)^{\mathrm{T}} \left( R\otimes I_{q} \right) \tilde{v}\notag\\ 
	&=\big(\operatorname{vec}\big(Q^{\mathrm{T}}\big)\big)^{\mathrm{T}} (R\otimes D^{\mathrm{T}}) \operatorname{vec}\left(\varPi\right) \notag\\
	&=\big(\operatorname{vec}\big(Q^{\mathrm{T}}\big)\big)^{\mathrm{T}}  \operatorname{vec}\big(D^{\mathrm{T}} \varPi R \big) \notag\\
	&=\operatorname{tr}\big(QD^{\mathrm{T}} \varPi R \big) = \operatorname{tr}\big(\varPi R QD^{\mathrm{T}} \big) \notag\\
	&=\operatorname{tr}\big(\varPi R \mathbf{1}_{N}d_{1}^{\mathrm{T}} \big) + \operatorname{tr}\big(\varPi R H \big),
	\end{align}
	where $ H = \sum_{i=2}^{N} \eta_{i}d_{i}^{\mathrm{T}} $. Observing that $ \big( r^{\mathrm{T}} \otimes I_{q}\big) v= \mathbf{0}  $, one thus concludes that $ \varPi R \mathbf{1}_{N}= \sum_{i=1}^{N} r_{i}\tilde{v}_{i} = \mathbf{0} $ by simple computation. Therefore, (\ref{Cx_v_1}) can be rewritten as follows,
	\begin{align}\label{cross terms}
	\left( C\tilde{x} \right)^{\mathrm{T}} \left( R\otimes I_{q} \right) \tilde{v} & = \operatorname{tr}\left(H \varPi R  \right)= \big(\operatorname{vec}\big(H^{\mathrm{T}}\big) \big)^{\mathrm{T}}\operatorname{vec}(\varPi R)   \notag\\
	& = \big(\operatorname{vec}\big(H^{\mathrm{T}}\big) \big)^{\mathrm{T}} \left( R\otimes I_{q} \right) \tilde{v}.
	\end{align}
	
	Note that $ \| \operatorname{vec}\big(H^{\mathrm{T}} \big) \|^{2} = \big(\operatorname{vec}\big(H^{\mathrm{T}}\big) \big)^{\mathrm{T}} \operatorname{vec}\big(H^{\mathrm{T}}\big) = \operatorname{tr}\big( H H^{\mathrm{T}} \big) = \| s \|^{2} $ holds. Then, by using equation (\ref{cross terms})  together with $ r_{i} < 1 $, $ i=1,2,\ldots,N $, we have the following inequality, 
	\begin{align}\label{inequality_1.2}
	-2\gamma_{2} \left( C\tilde{x} \right)^{\mathrm{T}} \left( R\otimes I_{q} \right) \tilde{v} & \leq \dfrac{\gamma_{2}^{2}}{\delta} \| \operatorname{vec}\big(H^{\mathrm{T}} \big) \|^{2} + \delta \| \tilde{v}\|^{2}  \notag\\
	& = \dfrac{\gamma_{2}^{2}}{\delta} \| s \|^{2} + \delta \| \tilde{v}\|^{2}.
	\end{align}
	Meanwhile, we have
	\begin{align} 		
	2 \left(C\tilde{x}  \right)^{\mathrm{T}} &\left( R\otimes I_{q} \right)CA_{c}(\tilde{x}-\Psi \tilde{\rho}) \notag\\ 
	& \leq \dfrac{\|C\|^{2}}{\delta} \|C\tilde{x}\|^{2} + \delta\|A_{c}\|^{2} \|\tilde{x}-\Psi \tilde{\rho}\|^{2}, \label{inequality_1.3} \\		
	\tilde{v}^{\mathrm{T}}  \big( R\otimes & I_{q} \big) CA_{c}(\tilde{x}-\Psi \tilde{\rho}) \notag\\ 
	& \leq \dfrac{\|C\|^{2}}{4\delta} \|\tilde{v}\|^{2} + \delta\|A_{c}\|^{2} \|\tilde{x}-\Psi \tilde{\rho}\|^{2}. \label{inequality_1.4}	
	\end{align}
	Substituting (\ref{inequality_1.1}), (\ref{inequality_1.2})--(\ref{inequality_1.4}) into (\ref{inequality_1}) leads to
	\begin{align} \label{dot_V1}
	\dot{V_{1}}\leq & -\Big( 2m-\dfrac{M^{2}+\|C\|^{2}}{\delta} \Big) \lVert C\tilde{x} \rVert^{2} \notag\\
	&- \Big(\gamma_{2} r_{\min}-\dfrac{5\delta}{4}- \dfrac{\|C\|^{2}}{4\delta} \Big) \|\tilde{v}\|^{2}   \notag\\
	& -\Big( \gamma_{1} \lambda_{2}-\dfrac{\gamma_{2}^{2}}{\delta} \Big) \|s\|^{2}  + 2 \delta\|A_{c}\|^{2} \|\tilde{x}-\Psi \tilde{\rho}\|^{2} . 
	\end{align}
	
	Since $ A_{c} $ is Hurwitz, there exists a positive definite matrix $ P_{c} $ such that $ P_{c}A_{c}+ A_{c}^{\mathrm{T}}P_{c} \leq -I $. Choose another positive definite function $ V_{2}= \left(\tilde{x}-\Psi \tilde{\rho} \right)^{\mathrm{T}} P_{c} \left(\tilde{x}-\Psi \tilde{\rho} \right) $. It follows from (\ref{dot_tilde_x}) that the dynamics of $ \tilde{x}-\Psi \tilde{\rho} $ satisfies 
	\begin{equation}
	\dot{\tilde{x}}-\Psi \dot{\tilde{\rho}} = A_{c}\left(\tilde{x}-\Psi \tilde{\rho} \right).
	\end{equation} 
	Thus, the derivative of $ V_{2} $ along the trajectory of (\ref{change system}) is inferred as
	\begin{align} \label{dot_V2}
	\dot{V_{2}} & = \left(\tilde{x}-\Psi \tilde{\rho} \right)^{\mathrm{T}} \big( P_{c}A_{c}+ A_{c}^{\mathrm{T}}P_{c}\big)  \left(\tilde{x}-\Psi \tilde{\rho} \right) \notag \\
	& \leq -\|\tilde{x}-\Psi \tilde{\rho}\|^{2}.
	\end{align}
	
	Consider the Lyapunov function candidate $ V(\tilde{x},\tilde{\rho},\tilde{v})= V_{1} + \delta_{c} V_{2} $ where $ \delta_{c}= 2\delta \|A_{c}\|^{2} +1 $. 
	Then by combining (\ref{dot_V1}) and (\ref{dot_V2}), the derivative of $ V $ along the trojectory of (\ref{change system}) satisfies,
	\begin{align} \label{inequality_3}
	\dot{V}\leq & -\Big( 2m-\dfrac{M^{2}+\|C\|^{2}}{\delta} \Big) \lVert C\tilde{x} \rVert^{2} \notag\\
	&- \Big(\gamma_{2} r_{\min}-\dfrac{5\delta}{4}- \dfrac{\|C\|^{2}}{4\delta} \Big) \|\tilde{v}\|^{2}  \notag \\
	& -\Big( \gamma_{1} \lambda_{2}-\dfrac{\gamma_{2}^{2}}{\delta} \Big) \|s\|^{2} - \|\tilde{x}-\Psi \tilde{\rho}\|^{2}.  
	\end{align}
	Choose the coupling gains $ \gamma_{1} $ and $ \gamma_{2} $ such that the following inequalities are satisfied,
	\begin{equation}\label{gains condition}
	\left\{\begin{array}{l}
	{2m-\dfrac{M^{2}+\|C\|^{2}}{\delta} > 0,} \\
	{\gamma_{2} r_{\min}-\dfrac{5\delta}{4}- \dfrac{\|C\|^{2}}{4\delta} > 0,} \\
	{\gamma_{1} \lambda_{2}-\dfrac{\gamma_{2}^{2}}{\delta}  > 0}.
	\end{array}\right.
	\end{equation}
	One then has
	\begin{equation}\label{inequality_final}
	\dot{V} \leq -\varepsilon\left(\lVert C\tilde{x} \rVert^{2} + \|\tilde{v}\|^{2} + \|\tilde{x}-\Psi \tilde{\rho}\|^{2} \right), 
	\end{equation}
	where $ \varepsilon=\min\left\lbrace 2m-\frac{M^{2}+\|C\|^{2}}{\delta}, \gamma_{2} r_{\min}-\frac{5\delta}{4} - \frac{\|C\|^{2}}{4\delta},1 \right\rbrace $. In fact, by selecting appropriate parameters $ \delta > \frac{M^{2}+\|C\|^{2}}{2m} $, we then can successively choose sufficiently large $ \gamma_{2}> \frac{5\delta^{2}+\|C\|^{2}}{4\delta r_{\min}} $ and $ \gamma_{1} > \frac{\gamma_{2}^{2}}{\lambda_{2}\delta} $ to ensure the inequalities in (\ref{gains condition}). In other words, (\ref{gains condition}) can always be satisfied by choosing large enough constants $ \gamma_{1} $ and $ \gamma_{2} $.
	
	Note that $ V(\tilde{x},\tilde{\rho},\tilde{v}) $ can be rewritten as 
	\begin{equation*}
	V = \left( \begin{array}{ccc} {C\tilde{x}}  \\ {\tilde{v}} \\ {\tilde{x}-\Psi \tilde{\rho}}  \end{array}\right)^{\mathrm{T}} F \left( \begin{array}{ccc} {C\tilde{x}}  \\ {\tilde{v}} \\ {\tilde{x}-\Psi \tilde{\rho}}  \end{array}\right) 
	\end{equation*}
	by defining the following positive definite matrix 
	\begin{equation*}
	F=\left(\begin{array}{ccc}
	{\left(R\otimes I_{q} \right) } & {\frac{1}{2}\left(R\otimes I_{q} \right)} & {\mathbf{0}} \\
	{\frac{1}{2}\left(R\otimes I_{q} \right)} & {\frac{1}{2}\left(R\otimes I_{q} \right)} & {\mathbf{0}}\\
	{\mathbf{0}} & {\mathbf{0}} & {\delta_{c} P_{c}}
	\end{array}\right).
	\end{equation*}
	Let $ \mu $ denote the maximum eigenvalue of $ F $. One then has $ V\leq \mu \left(\lVert C\tilde{x} \rVert^{2} + \|\tilde{v}\|^{2} + \|\tilde{x}-\Psi \tilde{\rho}\|^{2} \right) $ from the definition of $ V(\tilde{x},\tilde{\rho},\tilde{v}) $. Thus, (\ref{inequality_final}) can be rewritten as $ \dot{V}\leq -\frac{\varepsilon}{\mu}V $. With this fact, one can claim that $ \lim_{t\to\infty}C\tilde{x} = \mathbf{0} $, $ \lim_{t\to\infty}\tilde{v} = \mathbf{0} $ and $ \lim_{t\to\infty}\left( \tilde{x}- \Psi \tilde{\rho}\right) = \mathbf{0} $ exponentially, with a convergence rate no less than $ \varepsilon/\mu $. 
	
	Note that $ \tilde{\rho} = C\tilde{x}- C\left(\tilde{x}-\Psi \tilde{\rho} \right) $. By referring to $ \lim_{t\to\infty} C\tilde{x}= \mathbf{0} $ and $ \lim_{t\to\infty}\left( \tilde{x}- \Psi \tilde{\rho}\right) = \mathbf{0} $ exponentially, one then obtains that $ \lim_{t\to\infty}\tilde{\rho} = \mathbf{0} $ exponentially. Thus, using $ \lim_{t\to\infty} \left(\tilde{x}-\Psi \tilde{\rho} \right) = \mathbf{0} $, one can obtain that $ \lim_{t\to\infty}\tilde{x} = \mathbf{0} $ exponentially. Therefore, the exponential stability of system (\ref{subsystem}) is established.
	
	Step 2: The exponential stability of the following system (\ref{perturbed system}) is presented, 
	\begin{equation}\label{perturbed system}
	\dot{\xi}=g(\xi)+\kappa(t,\xi).
	\end{equation}
	Note that (\ref{perturbed system}) can be interpreted as the perturbed system of (\ref{subsystem}), where the perturbation term $ \kappa(t,\xi) $ satisfies $ \kappa(t,\bar{\xi})= \mathbf{0} $. Moreover, $ \kappa(t,\xi) \leq \sigma(t) \|\xi-\bar{\xi}\| $ with $ \sigma(t)= M\|C\| \max_{i}| r_{i}^{-1}-(z_{i}^{i}(t))^{-1} | $. Since it is proved that $ \lim_{t\to\infty} Z_{N}^{-1}(t) = R^{-1} $ exponentially, we have $ \lim_{t\to\infty} \sigma(t)= 0 $ exponentially. Then it follows from Lemma \ref{lemma 4} that the equilibrium point $ \bar{\xi} $ of the perturbed system (\ref{perturbed system}) is exponentially stable. 
	
	 Step 3: The exponential stability of system (\ref{perturbed form}) at the equilibrium point $ \bar{\xi}=\operatorname{col}(\bar{x},\bar{\rho},\bar{v}) $ is established. Since $ \nabla\tilde{f} $ is globally Lipschitz by Assumption \ref{cost function assumption}, we learn that $ G(\xi,\omega)= g(\xi)+\kappa(t,\xi) + \omega(t) $ is globally Lipschitz in $ (\xi, \omega) $. The boundness of $ \nabla\tilde{f}(\bar{Y}) $ suggests that $ \omega(t) $ is bounded. Then it follows from Lemma 4.6 in \cite{khalil2002nonlinear} that the system (\ref{perturbed form}) is input-to-state stable (ISS), i.e., for any initial state $ \xi(t_{0}) $ and any bounded input $ \omega(t) $, there exist a class $ \mathcal{KL} $ function $ \varphi $ and a class $ \mathcal{K} $ function $ \psi $ such that the solution $ \xi(t) $ satisfies $ \|\xi-\bar{\xi}\| \leq \varphi(\|\xi(t_{0})-\bar{\xi}\|,t-t_{0}) + \psi (\sup_{t_{0}\leq \tau \leq t} \omega(\tau)) $ for all $ t \geq t_{0} $. 
	Recalling the fact that $ \lim_{t\to\infty} \big(Z_{N}^{-1}(t) - R^{-1} \big) = \mathbf{0} $ exponentially, we have $ \lim_{t\to\infty} \omega(t) = \mathbf{0} $ exponentially. Thus it can be shown that $ \xi $ exponentially converges to $ \bar{\xi} $ by the property of ISS given in \cite{khalil2002nonlinear}. Therefore, we further claim that $ Y $ exponentially converges to $ \bar{Y}= \mathbf{1}_{N}\otimes y^{*} $, with $ y^{*} $ being the solution of the global cost function. The proof is thus completed.
	\end{proof}			
\end{theorem}

\begin{remark}
	Inspired by work \cite{zhu2018continuous}, we deal with the difficulties generated by the asymmetric information transmission caused by the weight-unbalanced directed graphs in virtue of the useful results from Kronecker matrix algebra and direct sum operation of vectors instead of the commonly used orthogonal transformation. 
\end{remark}

\begin{remark}
	It can be seen from the inequalities in (\ref{gains condition}) that the choice of control gains $ \gamma_{1} $ and $ \gamma_{2} $ only depends on the minimum value of the elements in the left eigenvector, instead of the exact value of the left eigenvector as in \cite{li2017distributed}. Therefore, once we can obtain the lower bound of the minimum value without any global information, it can be directly applied to the controller design in this work. Furthermore, the inequalities in (\ref{gains condition}) can always be guaranteed as long as control parameters $ \gamma_{1} $ and $ \gamma_{2} $ are chosen to be large enough. 
\end{remark}


\subsection{Distributed observer-based output feedback controller}

It should be noted that state measurements may be unavailable in practical scenarios. In this subsection, the previous distributed state feedback controller is extended to a distributed observer-based output feedback controller. More specifically, the newly proposed output feedback control law is given as follows,
\begin{equation}\label{observer based control law}	
\begin{aligned}	
u_{i} & =-K_{i}\hat{x}_{i}+\Upsilon_{i}\omega_{i}-\left(\Phi_{i}-K_{i}\Psi_{i} \right) \rho_{i},\\
\dot{\hat{x}}_{i} & = A_{i}\hat{x}_{i} + B_{i}u_{i} + H_{i}(y_{i}-C_{i}\hat{x}_{i}), \\
\dot{\rho}_{i} & =\omega_{i} :=-\dfrac{\nabla f_{i}(y_{i})}{z_{i}^{i}} -\gamma_{1}\sum_{j=1}^{N}a_{i j}\left(y_{i}-y_{j} \right)-\gamma_{2} v_{i}, \\
\dot{v}_{i} & =\gamma_{1} \sum_{j=1}^{N} a_{i j}\left(y_{i}-y_{j} \right), \\
\dot{z}_{i} & =-\sum_{j=1}^{N} a_{i j}\left(z_{i}-z_{j} \right) ,
\end{aligned}	
\end{equation}
where $ \hat{x}_{i} $ is the estimation of state $ x_{i} $, $ H_{i}\in\mathbb{R}^{n_{i}\times q} $ is an observer feedback matrix such that $ A_{i}-H_{i}C_{i} $ is Hurwitz, and the remaining variables are defined to be the same as those in control law (\ref{control law}).

Substituting the above control law into system dynamics (\ref{dynamics}) yields the following compact form of the closed-loop system,
\begin{equation}	\label{observer based closed-loop system}
\begin{aligned}	
\dot{x} &\! =Ax-BK\hat{x} + B\Upsilon \dot{\rho}-B\left( \Phi-K\Psi \right) \rho ,\\
\dot{\hat{x}} &\! =\!\left( A\!-\!BK\right) \hat{x} \!+\! B\Upsilon \dot{\rho}\!-\!B\left( \Phi\!-\!K\Psi \right) \rho \!+\! H\left(Y\!-\!C\hat{x} \right) ,\\
\dot{\rho} & \!=-\left( Z_{N}^{-1}\otimes I_{q} \right)\nabla\tilde{f}(Y)-\gamma_{1} \left( \mathcal{L}\otimes I_{q}\right)Y-\gamma_{2} v  , \\
\dot{v} & \!=\gamma_{1} \left( \mathcal{L}\otimes I_{q}\right)Y,  \\
\dot{z} & \!=-\left(\mathcal{L}\otimes I_{N} \right)z,  
\end{aligned}	
\end{equation}
where $ \hat{x}=\mathrm{col}(\hat{x}_{1},\hat{x}_{2},\ldots,\hat{x}_{N}) $,  $ H=\mathrm{diag}(H_{1},H_{2},\ldots,H_{N}) $, and the remaining terms are defined as those in the closed-loop system (\ref{closed-loop system}). Define the estimation error variable $ \breve{x}= x-\hat{x} $. Then (\ref{observer based closed-loop system}) can be rewritten as follows,
\begin{equation}
	\begin{aligned}	\label{alternative observer based closed-loop system}
	\dot{x} & =Ax-BK\left(x-\breve{x} \right)  + B\Upsilon \dot{\rho}-B\left( \Phi-K\Psi \right) \rho ,\notag\\
	\dot{\breve{x}} & =\left( A-HC\right) \breve{x}  ,\notag\\
	\dot{\rho} & =-\left( Z_{N}^{-1}\otimes I_{q} \right)\nabla\tilde{f}(Y)-\gamma_{1} \left( \mathcal{L}\otimes I_{q}\right)Y-\gamma_{2} v  , \\
	\dot{v} & =\gamma_{1} \left( \mathcal{L}\otimes I_{q}\right)Y, \notag \\
	\dot{z} & =-\left(\mathcal{L}\otimes I_{N} \right)z.\notag  
	\end{aligned}
\end{equation}

The second main result of this paper is presented below.
\begin{theorem}\label{theorem 2}
	Consider system (\ref{dynamics}). Suppose Assumptions \ref{cost function assumption}-\ref{matrix equations assumption} hold and $ \left(A_{i},C_{i} \right) $ is detectable. Then there exist two positive constants  $ \gamma_{1}$ and $\gamma_{2} $ such that the output $ Y $ exponentially reaches the optimal value $ \bar{Y}= \mathbf{1}_{N}\otimes y^{*} $ of problem (\ref{equivalent problem}) under the output feedback controller  (\ref{observer based control law}), with $ y^{*} $ be the optimal solution to problem (\ref{objetive}).
	
	\begin{proof}
		The proof is given in the Appendix.
	\end{proof}
\end{theorem}


\section{Illustrative Examples} \label{section simulation results}
In this section, the effectiveness of the proposed control laws is illustrated by two illustrative examples. We will start with a simplified but practical example in RLC networks. 

\begin{figure}[!t] 
	\centering  
	\includegraphics[height=4.3cm, width=8.7cm]{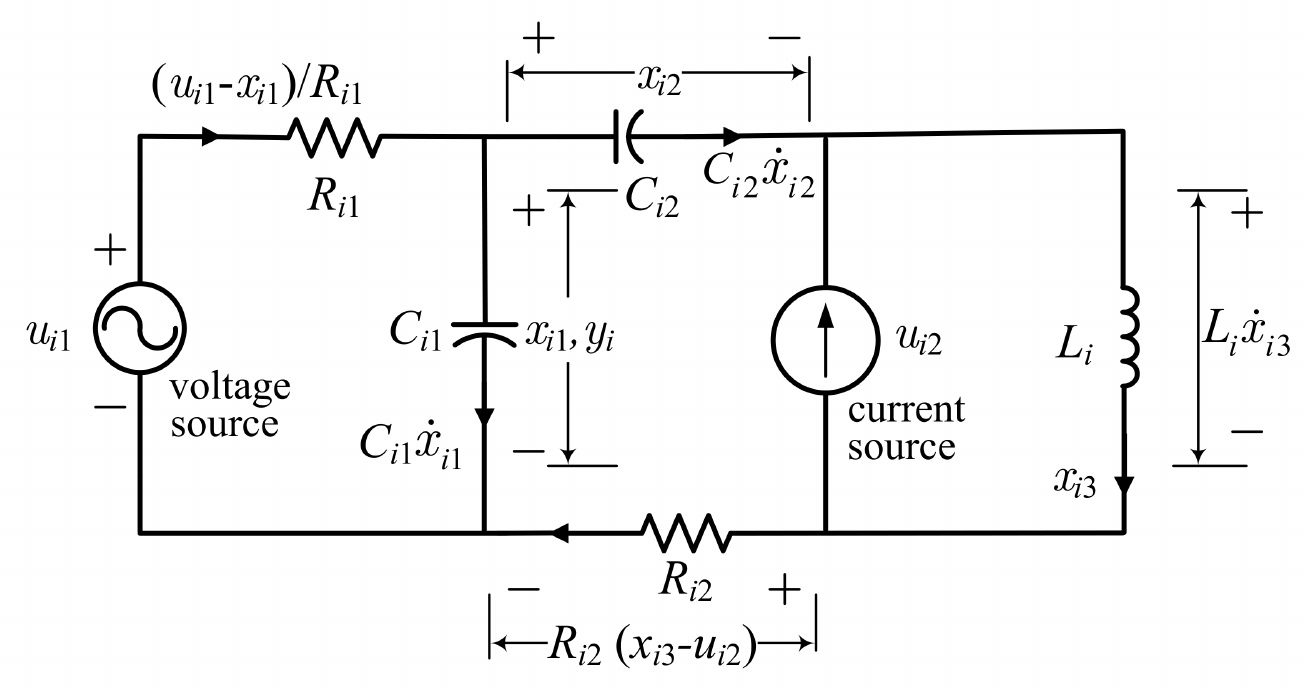} 
	\caption{RLC network.}
	\label{Fig_RLC} 
\end{figure}

\textbf{Example 1.}
Consider the RLC network depicted in Fig. \ref{Fig_RLC}, which is a modification of Figure 2.16 in \cite{chen1999linear}. It consists of the voltage source $ u_{i1} $, current source $ u_{i2} $, two resistors $ R_{i1} $ and $ R_{i2} $, inductor $ L_{i} $, and two capacitors $ C_{i1} $ and $ C_{i2} $. The capacitor voltages and the inductor current will be assigned as state variables $ x_{i1} $, $ x_{i2} $ and $ x_{i3} $, respectively. Then we can apply Kirchhoff's current and voltage laws to establish the following equations,
\begin{equation*}
\begin{aligned}
x_{i3} & = \frac{u_{i1}-x_{i1}}{R_{i1}} - C_{i1}\dot{x}_{i1} + u_{i2}, \\
x_{i3} & = C_{i2}\dot{x}_{i2} + u_{i2}, \\
x_{i1} & = x_{i2} + L_{i}\dot{x}_{i3} + R_{i2}(x_{i3}-u_{i2}), \\
y_{i}  & = x_{i1}.
\end{aligned}
\end{equation*}
By defining $ x_{i}=\operatorname{col}(x_{i1}, x_{i2}, x_{i3}) $ and $ u_{i}=\operatorname{col}(u_{i1}, u_{i2}) $, the state-space description of the RLC network
takes the linear dynamics form (\ref{dynamics}) with
\begin{equation}\label{RLC_matrix}
\begin{aligned}
A_{i} & = \left[\begin{array}{ccc}
{-\dfrac{1}{C_{i1}R_{i1}}} & {0} & {-\dfrac{1}{C_{i1}}} \\
{0} & {0} & {\dfrac{1}{C_{i2}}} \\
{\dfrac{1}{L_{i}}} & {-\dfrac{1}{L_{i}}} & {-\dfrac{R_{i2}}{L_{i}}} 
\end{array}\right], \\
B_{i} & = \left[\begin{array}{cc}
{\dfrac{1}{C_{i1}R_{i1}}} & {\dfrac{1}{C_{i1}}} \\
{0} & {-\dfrac{1}{C_{i2}}} \\
{0} & {\dfrac{1}{L_{i}}} 
\end{array}\right], 
C_{i}  = \left[\begin{array}{ccc}
{1} & {0} & {0}
\end{array}\right].
\end{aligned}
\end{equation}

\begin{figure}[t] 
	\centering 
	\begin{tikzpicture}[> = stealth, 
	shorten > = 1pt, 
	auto,
	node distance = 3cm, 
	semithick 
	,scale=0.6,auto=left,every node/.style={circle,fill=blue!5,draw=black!80,text centered}]
	\centering
	\node (n1) at (0,0)		{1};
	\node (n2) at (2,0)  	{2};
	\node (n3) at (4,0) 	{3};
	\node (n4) at (6,0) 	{4};

	\draw[->,black!80] (n3) to [out=135,in=45] (n1);
	\draw[->,black!80] (n1)-- (n2);
	\draw[->,black!80] (n2)-- (n3);
	\draw[->,black!80] (n3)-- (n4);
	\draw[->,black!80] (n4) to [out=-135,in=-45] (n2);
	
	\end{tikzpicture} 
	\caption{Communication between RLC networks.} 
	\label{Fig_topology1}
\end{figure}
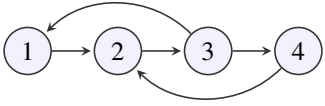

Consider a multi-agent system of 4 agents described by the RLC network as shown in Fig. \ref{Fig_RLC}, the communication topology among agents is given in Fig. \ref{Fig_topology1}. Assume that local cost functions of RLC circuits are described by the following quadratic form functions,
\begin{equation*}
\begin{aligned}
f_{1} & =0.2y^{2}-2y+1, \quad f_{2}=0.4y^{2}+y+2, \\
f_{3} & =0.6y^{2}-3y-1, \quad f_{4}=0.8y^{2}+y+1.
\end{aligned}
\end{equation*}

\begin{table}[!htbp]
	\centering
	\caption{Parameters of RLC networks}\label{table 1}
	\begin{tabular}{cccccc}
		\hline \hline 
		i & $ R_{i1} $ & $ R_{i2} $ & $ L_{i} $ & $ C_{i1} $ & $ C_{i2} $ \\
		\hline 1 & 2 & 1 & 3 & 1 & 2 \\
		2 & 1 & 2 & 2 & 3 & 1 \\
		3 & 0.5 & 2 & 1 & 0.5 & 3 \\
		4 & 3 & 0.5 & 2 & 1 & 0.5\\
		\hline \hline
	\end{tabular}
\end{table}

Then, it is calculated that the optimal value of the global cost function is $ y^{*}=1.5 $. Moreover, it can be verified that Assumptions \ref{cost function assumption} and \ref{graph assumption} are satisfied. Furthermore, to ensure the solvability of linear matrix equations (\ref{matrix equations}) corresponding to dynamics matrices (\ref{RLC_matrix}), for each agent $ i, i=1,2,3,4 $, the parameters of (\ref{RLC_matrix}) are given in TABLE \ref{table 1} so that Assumption \ref{matrix equations assumption} is satisfied. Then by Lemma \ref{lemma_matrix equations}, the solution triplets $ \left( \Upsilon_{i}, \Phi_{i}, \Psi_{i} \right), i=1,2,3,4 $ of linear matrix equations (\ref{matrix equations}) can be chosen as follows,
\begin{align*}\tiny
\Upsilon_{1} & = \Big[\begin{array}{cc} {8}  \\ {-1}  \end{array}\Big], \quad \Upsilon_{2}= \Big[\begin{array}{cc} {0.5}  \\ {1.5}  \end{array}\Big], \quad \Upsilon_{3}= \Big[\begin{array}{cc} {0.8}  \\ {-0.6}  \end{array}\Big], \\
\Upsilon_{4} & = \Big[\begin{array}{cc} {7.5}  \\ {-0.5}  \end{array}\Big], \quad \Phi_{1}= \Big[\begin{array}{cc} {-1}  \\ {0.5}  \end{array}\Big], \quad 
\Phi_{2}= \Big[\begin{array}{cc} {-1}  \\ {-1.5}  \end{array}\Big],\\
\Phi_{3} & = \Big[\begin{array}{cc} {-1}  \\ {0.2}  \end{array}\Big], \quad \Phi_{4}=\Big[\begin{array}{cc} {-1}  \\ {1}  \end{array}\Big], \quad 
\Psi_{1} =\Bigg[\begin{array}{ccc} {1}  \\ {1} \\ {-0.5} \end{array}\Bigg],\\
\Psi_{2} & =\Bigg[\begin{array}{ccc} {1}  \\ {-0.5} \\ {1.5}  \end{array}\Bigg], \quad \Psi_{3} =\Bigg[\begin{array}{ccc} {1}  \\ {1.2} \\ {-0.2} \end{array}\Bigg], \quad 
\Psi_{4} =\Bigg[\begin{array}{ccc} {1}  \\ {0.5} \\ {-1} \end{array}\Bigg].
\end{align*}
Finally, the matrices $ K_{i} $ are chosen as follows such that $ A_{i}-B_{i}K_{i}, i=1,2,3,4 $ are Hurwitz,
\begin{align*}\tiny
K_{1} & =\Big[\begin{array}{ccc} {1} & {2} & {-2}  \\ {-1} & {0} & {1}  \end{array}\Big], \quad
K_{2} =\Big[\begin{array}{ccc} {0.5} & {2} & {-1}  \\ {-2} & {0} & {2}  \end{array}\Big],\\
K_{3} & =\Big[\begin{array}{ccc} {2} & {-1} & {-2}  \\ {0} & {-3} & {3}  \end{array}\Big], \quad
K_{4} =\Big[\begin{array}{ccc} {-2} & {1} & {2}  \\ {0} & {-1} & {2}  \end{array}\Big].
\end{align*}

\begin{figure}[t] 
	\begin{center}
		\includegraphics[height=3.2cm, width=8.5cm]{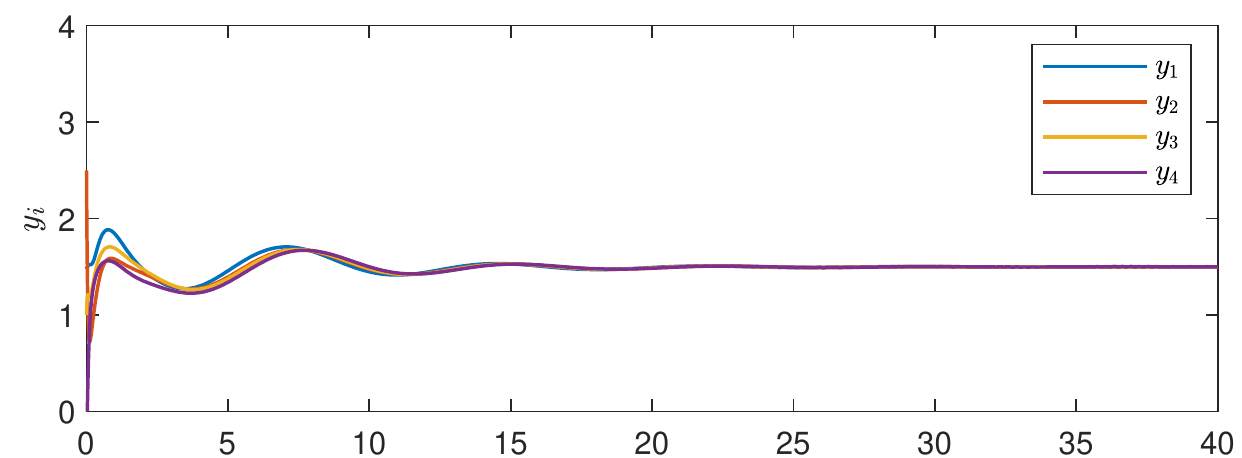} 
		\caption{Convergence performance of RLC networks.}
		\label{Fig_RLC_simu} 
	\end{center}
\end{figure}

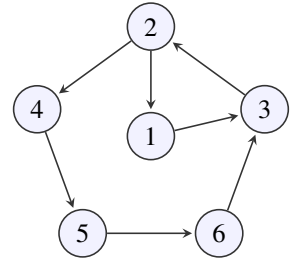
\begin{figure}[t] 
	\centering 
	\begin{tikzpicture}[> = stealth, 
	shorten > = 1pt, 
	auto,
	node distance = 3cm, 
	semithick 
	,scale=0.55,auto=left,every node/.style={circle,fill=blue!5,draw=black!80,text centered}]
	\centering
	\node (n2) at (0,0)		{2};
	\node (n4) at (-2.75,-2)  	{4};
	\node (n3) at (2.75,-2) 	{3};
	\node (n5) at (-1.65,-5) 	{5};
	\node (n6) at (1.65,-5) 	{6};
	\node (n1) at (0,-2.65) 	{1};
	\draw[->,black!80] (n1)--(n3);
	\draw[->,black!80] (n2)--(n1);
	\draw[->,black!80] (n2)--(n4);
	\draw[->,black!80] (n3)--(n2);
	\draw[->,black!80] (n4)--(n5);
	\draw[->,black!80] (n5)--(n6);
	\draw[->,black!80] (n6)--(n3);
	\end{tikzpicture} 
	\caption{Weight-unbalanced directed network.} 
	\label{Fig_topology}
\end{figure}

The initial value $ v(0) $ is chosen as $ v(0)= \mathbf{0} $. The other initial conditions of the closed-loop system are randomly chosen in the closed interval $ [-4,6] $. The simulation result is shown in Fig. \ref{Fig_RLC_simu}. It can be observed that the trajectories of all outputs $ y_{i} $'s converge to the global minimizer $ y^{*}=1.5 $. Thus the closed-loop system composed of (\ref{dynamics}) and (\ref{control law}) achieves optimal output consensus eventually. \qed\\

Next, we will provide another example to compare the convergence performance under state feedback control law (\ref{control law}) for different control gains and illustrate the effectiveness of observer-based output feedback control law (\ref{observer based control law}).

\textbf{Example 2.}
Consider a group of 6 agents with the unbalanced directed network $ \mathcal{G} $ depicted in Fig. \ref{Fig_topology}. For agents $ i=1,2,\ldots,6 $, the local cost functions are respectively given as follows,
\begin{equation*}
	\begin{aligned}
		f_{1} & =\sin(0.2y-(\pi/2)), \quad f_{2}=0.2\cos(\ln(y^{2}+4)-0.2), \\
		f_{3} & =0.1(y+0.3)^{2}+0.2(y-2)^{2}, \quad f_{4}=0.4y^{2}\ln(5+y^{2}),\\
		f_{5} & =0.2y^{2}(\ln(y^{2}+1)+1), \quad f_{6}=0.3y^{2}/\sqrt{y^{2}+5}.
	\end{aligned}
\end{equation*}
It can be verified that Assumptions \ref{cost function assumption} and \ref{graph assumption} are satisfied. Therefore, the strong convexity of the global cost function guarantees that the global minimizer $ y^{*}=0.286 $ is unique.

The dynamics of agents are described by (\ref{dynamics}) with
\begin{align*}
	A_{1} & = A_{2} = \Big[\begin{array}{cc} {0} & {1}  \\ {0} & {0}  \end{array}\Big], \quad B_{1}=B_{2}= \Big[\begin{array}{cc} {0} & {1}  \\ {1} & {-2}  \end{array}\Big], \\
	C_{1} & = C_{2}= \big[\begin{array}{c} {1} ~~ {1}  \end{array}\big], \quad A_{3}=A_{4}= \Big[\begin{array}{cc} {0} & {-1}  \\ {1} & {-2}  \end{array}\Big],\\
	B_{3} & = B_{4}= \Big[\begin{array}{cc} {1} & {0}  \\ {3} & {-1}  \end{array}\Big], \quad C_{3}=C_{4}= \big[\begin{array}{c} {-1} ~~ {1}  \end{array}\big],\\
	A_{5} & = A_{6}= \Bigg [\begin{array}{ccc} {0} & {1} & {0} \\ {0} & {0} & {1} \\ {0.5} & {1} & {-2} \end{array}\Bigg], \quad B_{5}=B_{6}= \Bigg[\begin{array}{ccc} {1} & {0}  \\ {0} & {1} \\ {1} & {0} \end{array} \Bigg],\\
	C_{5} & = C_{6}= \left[\begin{array}{c} {1} ~~ {-1} ~~ {1} \end{array}\right].
\end{align*}

Note that Assumption \ref{matrix equations assumption} is also satisfied. Then by Lemma \ref{lemma_matrix equations}, the solution triplets $ \left( \Upsilon_{i}, \Phi_{i}, \Psi_{i} \right), i=1,2,\ldots,6 $ of linear matrix equations (\ref{matrix equations}) can be chosen as follows,
\begin{align*}
	\Upsilon_{1} & = \Upsilon_{2}=\big[\begin{array}{cc} {1.5}  & {0.5}  \end{array}\big]^{\operatorname{T}}, \quad \Upsilon_{3}=\Upsilon_{4}=\big[\begin{array}{cc} {-0.5}  & {-2}  \end{array}\big]^{\operatorname{T}}, \\
	\Upsilon_{5} & = \Upsilon_{6}=\big[\begin{array}{cc} {0}  & {-1}  \end{array}\big]^{\operatorname{T}}, \quad \Phi_{1}=\Phi_{2}=\big[\begin{array}{cc} {1}  & {0.5}  \end{array}\big]^{\operatorname{T}},\\
	\Phi_{3} & = \Phi_{4}=\big[\begin{array}{cc} {-0.5}  & {0}  \end{array}\big]^{\operatorname{T}}, \quad \Phi_{5}=\Phi_{6}=\big[\begin{array}{cc} {-1}  & {0}  \end{array}\big]^{\operatorname{T}},\\
	\Psi_{1} & = \Psi_{2}=\big[\begin{array}{cc} {0.5}  & {0.5}  \end{array}\big]^{\operatorname{T}}, \quad \Psi_{3}=\Psi_{4}=\big[\begin{array}{cc} {-0.5}  & {0.5}  \end{array}\big]^{\operatorname{T}},\\
	\Psi_{5} & = \Psi_{6}=\big[\begin{array}{ccc} {0}  & {-1} & {0}  \end{array}\big]^{\operatorname{T}}.
\end{align*}

Furthermore, the matrices $ K_{i} $ and $ H_{i} $ are respectively chosen as follows such that $ A_{i}-B_{i}K_{i} $ and $ A_{i}-H_{i}C_{i}, i=1,2,\ldots,6 $ are Hurwitz,
\begin{align*}
K_{1} & = K_{2}=\Big[\begin{array}{cc} {3} & {5}  \\ {1.5} & {1}  \end{array}\Big], \quad K_{3}=K_{4}=\Big[\begin{array}{cc} {0.75} & {-1}  \\ {1.25} & {-4}  \end{array}\Big],\\
K_{5} & = K_{6}=\Big[\begin{array}{ccc} {2.167} & {1} & {0.333}  \\ {0} & {3} & {1}  \end{array}\Big], \quad H_{1} = H_{2}=\Big[\begin{array}{cc} {1}  \\ {2}  \end{array}\Big],\\
H_{3} & = H_{4}=\Big[\begin{array}{cc} {-2}  \\ {-1}  \end{array}\Big], \quad H_{5} =H_{6}=\big[\begin{array}{ccc} {4}  & {3} & {2}  \end{array}\big]^{\operatorname{T}}.
\end{align*}

\begin{figure}[!t] 
	\centering  
	\includegraphics[height=9.2cm, width=8.5cm]{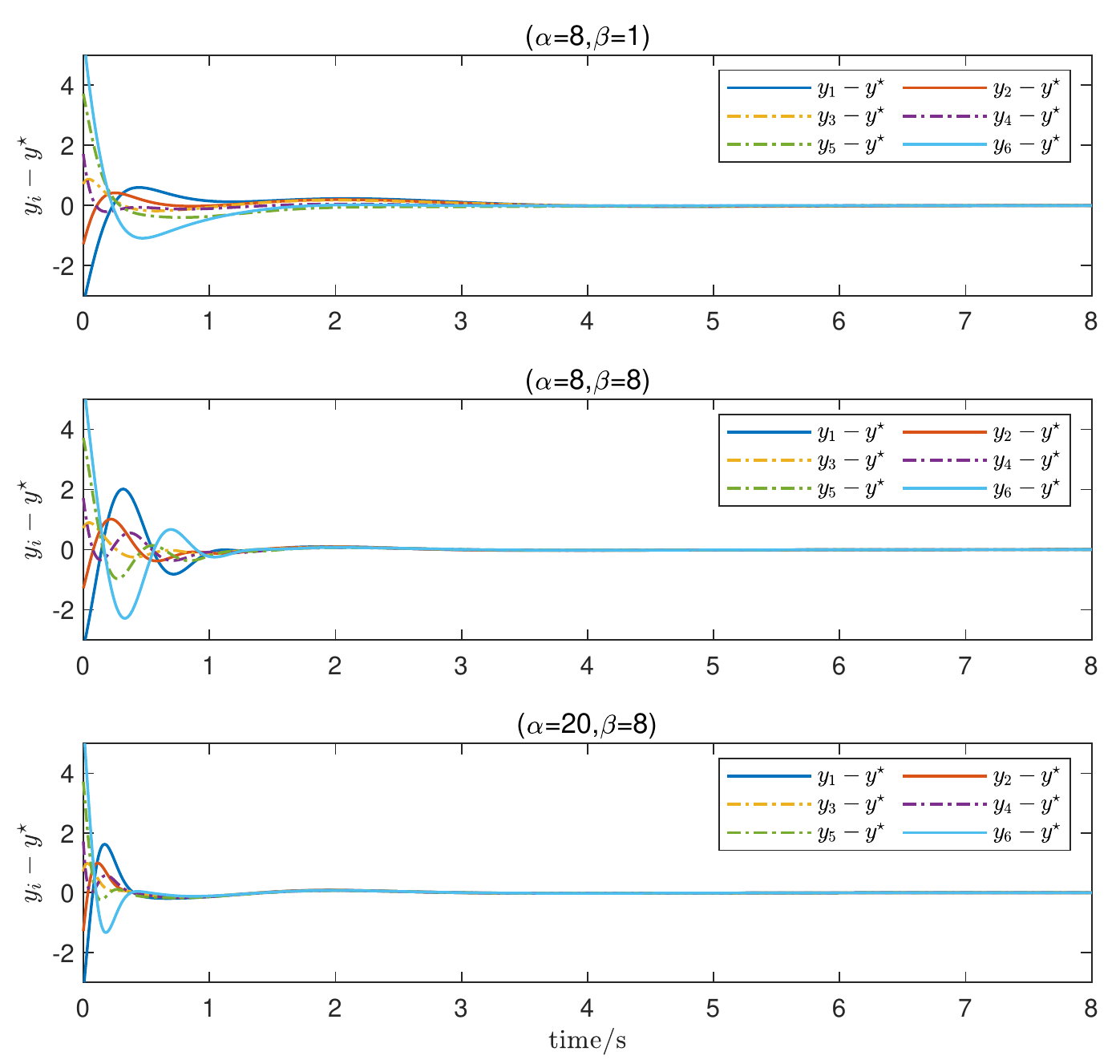} 
	\caption{Convergence performance comparisons under state feedback control law (\ref{control law}).}
	\label{Fig_statefeedback} 
\end{figure}
\begin{figure}[!t] 
	\centering 
	\includegraphics[height=6.7cm, width=8.5cm]{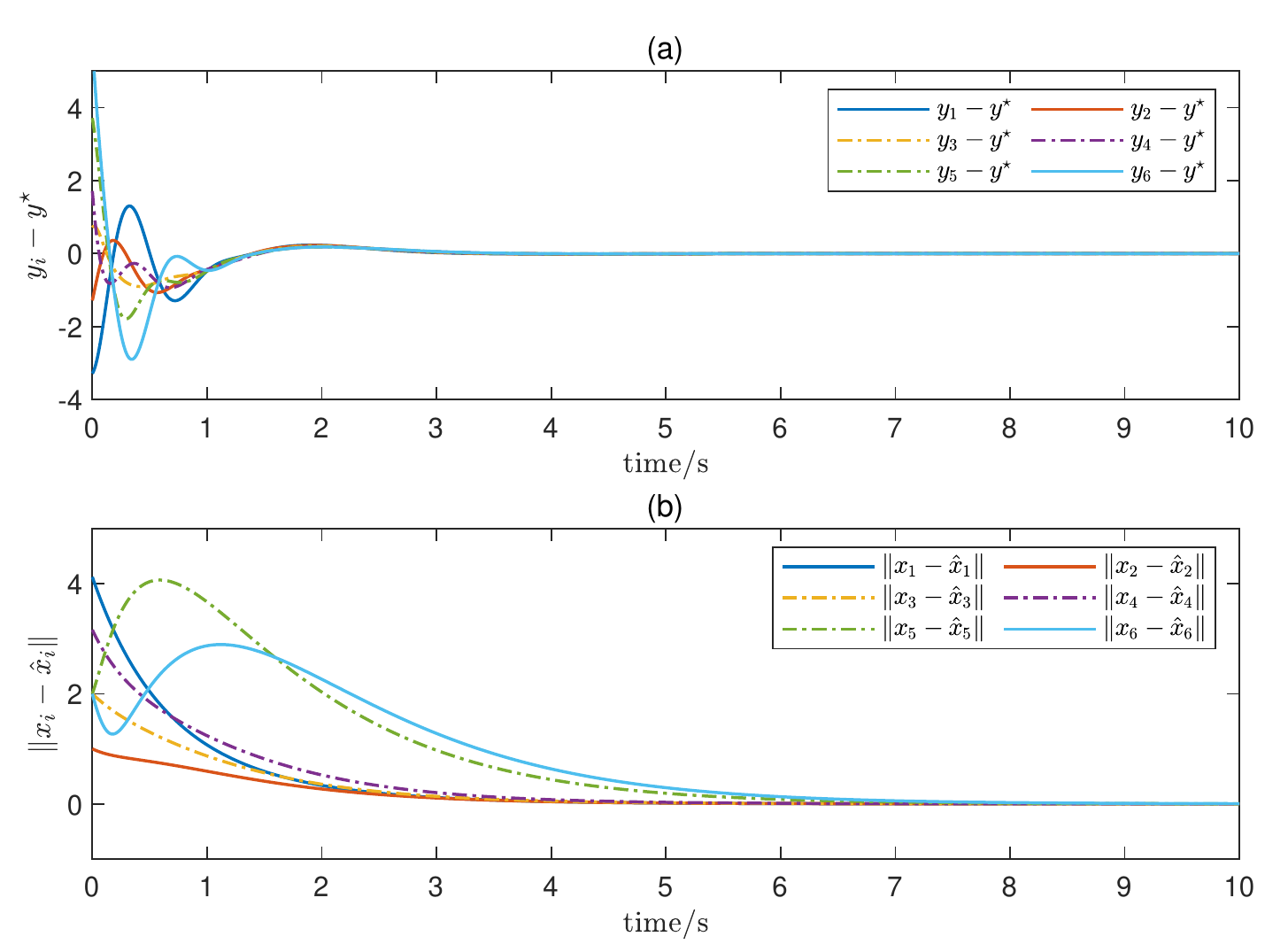} 
	\caption{Convergence performance under output feedback control law (\ref{observer based control law}).} 
	\label{Fig_outputfeedback}
\end{figure}
The initial values $ x(0), \hat{x}(0) $ and	$ \rho(0) $ can be arbitrarily chosen, while $ v(0) $ needs to satisfy $ v(0)= \mathbf{0} $. For convergence performance comparisons, three sets of control gains in (\ref{control law}) are given as $ (\gamma_{1}, \gamma_{2})= (8, 1), (\gamma_{1}, \gamma_{2})= (8, 8), $ and $ (\gamma_{1}, \gamma_{2})= (20, 8) $, respectively. The simulation results of the closed-loop control system via state feedback control law (\ref{control law}) are shown in Fig. \ref{Fig_statefeedback}. One can observe from the figure that the trajectories of outputs $ y_{i} $ converge to the global minimizer $ y^{*}=0.286 $. Moreover, by comparing the convergence performances in Fig. \ref{Fig_statefeedback}(a)-(c), it can be observered that the larger values of $ \gamma_{1} $ and $ \gamma_{2} $ will bring about faster convergence. 

To verify the effectiveness of the observer-based output feedback control law (\ref{observer based control law}), set $ (\gamma_{1}, \gamma_{2})= (8, 1) $ and the simulation results of the resulting closed-loop system are shown in Fig. \ref{Fig_outputfeedback}. Fig. \ref{Fig_outputfeedback}(a) shows that the outputs of all the agents would reach the optimal solution, while Fig. \ref{Fig_outputfeedback}(b) indicates that the observer states $ \hat{x}_{i} $ will eventually tend to $ x_{i}, i=1,2,\ldots,6 $, which is consistent with the theoretical result. \qed

\section{Conclusion} \label{section conclusion}
In this paper, we have studied the distributed optimal output consensus problem of heterogeneous linear multi-agent systems over weight-unbalanced directed networks. We have developed two novel distributed control laws. It is shown that the proposed control laws are able to ensure the agent outputs converge exponentially to the optimal solution under standard assumptions. Our results generalize the optimal output consensus problem of heterogeneous linear multi-agent systems from undirected networks to weight-unbalanced directed networks. Two illustrative examples have also been given to show the effectiveness of the proposed controllers. One of the possible future research topics is to further eliminate the dependence on global information of our current controller design. Another interesting topic is to extend the optimal output consensus problem over weight-unbalanced directed networks to the case that the agent dynamics are uncertain nonlinear systems.


\begin{appendix}[Proof of Theorem \ref{theorem 2}]\label{Appendix}

\section{Proof of Theorem \ref{theorem 2}} 
Define $ \zeta = \operatorname{col}(x,\breve{x},\rho,v) $, then the dynamics of $ \zeta $ is given as follows,
\begin{align} \label{observer based perturbed form}
\underbrace{
	\left(\begin{array}{c}
	{\dot{x}} \\
	{\dot{\breve{x}}}\\
	{\dot{\rho}}\\
	{\dot{v}}
	\end{array}\right)
}_{\dot{\zeta}}
&=\underbrace{
	\left(\begin{array}{c}
	{Ax-BK\left(x-\breve{x} \right)+B\Upsilon \dot{\rho}-B\left( \Phi-K\Psi \right) \rho}\\
	{\left( A-HC\right) \breve{x}}\\
	{-\left( R^{-1}\otimes I_{q} \right)\nabla\tilde{f}(Y)-\gamma_{1} \left( \mathcal{L}\otimes I_{q}\right)Y-\gamma_{2} v} \\
	{\gamma_{1} \left( \mathcal{L}\otimes I_{q}\right)Y}
	\end{array}\right)
}_{\bar{g}(\zeta)} \notag\\
&+\underbrace{
	\left(\begin{array}{c}
	{\mathbf{0}_{Nq}}\\
	{\mathbf{0}_{Nq}}\\
	{\left( \left( R^{-1}-Z_{N}^{-1}\right) \otimes I_{q} \right)\big( \nabla\tilde{f}(Y)-\nabla\tilde{f}(\bar{Y})\big) } \\
	{\mathbf{0}_{Nq}}
	\end{array}\right)
}_{\bar{\kappa}(t,\zeta)} \notag\\ 
&+\underbrace{
	\left(\begin{array}{c}
	{\mathbf{0}_{Nq}}\\
	{\mathbf{0}_{Nq}}\\
	{\left( \left( R^{-1}-Z_{N}^{-1}\right) \otimes I_{q} \right) \nabla\tilde{f}(\bar{Y}) } \\
	{\mathbf{0}_{Nq}}
	\end{array}\right)
}_{\bar{\omega}(t)} .
\end{align}

At first, we show that the output $ \bar{Y} $ at the equilibrium point $ \bar{\zeta}=\left(\bar{x},\bar{\breve{x}},\bar{\rho},\bar{v} \right)  $ of system $ \dot{\zeta}=\bar{g}(\zeta) $ is the solution of problem (\ref{equivalent problem}). Note that the point $ \bar{\zeta}=\left(\bar{x},\bar{\breve{x}},\bar{\rho},\bar{v} \right) $ satisfies
\begin{subequations}\label{output equilibrium point equations}	
	\begin{align}	
	\mathbf{0} & =\left( A-BK \right)\left(\bar{x}-\Psi \bar{\rho} \right) + BK\bar{\breve{x}} ,\label{output_equilibrium_x}\\
	\mathbf{0} & = \left( A-HC\right) \bar{\breve{x}}, \label{output_equilibrium_breve x} \\
	\mathbf{0} & =-\left( R^{-1}\otimes I_{q} \right)\nabla\tilde{f}(\bar{Y})-\gamma_{2}\bar{v}  , \label{output equilibrium_rho}\\
	\mathbf{0} & =\gamma_{1} \left( \mathcal{L}\otimes I_{q}\right)\bar{Y}.\label{output equilibrium_nu}  
	\end{align}	
\end{subequations}
Since $ A-HC $ is Hurwitz, we can obtain from (\ref{output_equilibrium_breve x}) that $ \bar{\breve{x}} = 0 $. Thus it can be inferred that the component $ \left(\bar{x},\bar{\rho},\bar{v} \right)  $ at the equilibrium point $ \bar{\zeta}=\left(\bar{x},0,\bar{\rho},\bar{v} \right)  $ of system $ \dot{\zeta}=\bar{g}(\zeta) $ coincides with that of system (\ref{perturbed form}), i.e., $ \left(\bar{x},\bar{\rho},\bar{v} \right)  $ satisfies equations (\ref{equilibrium point equations}). Then it can be shown that the output $ \bar{Y} $ at the equilibrium point of $ \dot{\zeta}=\bar{g}(\zeta) $ is the solution of problem (\ref{equivalent problem}) via a similar analysis as in Lemma \ref{Lemma equilibrium point}.

In what follows, the exponential stability of the following system is presented,
\begin{equation} \label{observer based perturbed system}
\dot{\zeta}=\bar{g}(\zeta). 
\end{equation}
To this end, transforming the equilibrium point of (\ref{observer based perturbed system}) to the origin by defining $ \tilde{x}=x-\bar{x} $, $ \tilde{\rho}=\rho-\bar{\rho} $ and $ \tilde{v}=v-\bar{v} $ leads to
\begin{equation}\label{observer based change system}	
\begin{aligned}	
\dot{\tilde{x}} & =A_{c} \left(\tilde{x}-\Psi \tilde{\rho} \right) +\Psi \dot{\tilde{\rho}} +BK\breve{x} , \\ 
\dot{\breve{x}} & = A_{o} \breve{x}, \\
\dot{\tilde{\rho}} & =-\left( R^{-1}\otimes I_{q} \right)h-\gamma_{1} \left(\mathcal{L}\otimes I_{q} \right)C\tilde{x}  -\gamma_{2}\tilde{v}  , \\
\dot{\tilde{v}} & =\gamma_{1} \left( \mathcal{L}\otimes I_{q}\right)C\tilde{x},  
\end{aligned}	
\end{equation}
where $ A_{c}=A-BK $, $ A_{o}=A-HC $ and $ h=\nabla\tilde{f}\left( C(\bar{x}+\tilde{x}) \right)- \nabla\tilde{f}\left( C\bar{x} \right) $.  

Reconsider the Lyapunov function candidate $ V(\tilde{x},\tilde{\rho},\tilde{v})= V_{1} + \delta_{c} V_{2} $, where $ V_{1} $, $ V_{2} $ and $ \delta_{c} $ are the same as those defined in the proof of Theorem \ref{Theorem 1}. The derivative of $ V $ along the trajectory of (\ref{observer based change system}) is given as follows,
\begin{equation*}
\begin{aligned}
\dot{V}= & -2\left( C\tilde{x} \right)^{\mathrm{T}}h- \gamma_{1} \left( C\tilde{x} \right)^{\mathrm{T}} \left( R\mathcal{L}\otimes I_{q} \right) \left( C\tilde{x} \right)  \\
& + 2 \left( C\tilde{x} \right)^{\mathrm{T}} \left( R\otimes I_{q} \right) CA_{c}\left( \tilde{x}-\Psi\tilde{\rho} \right) -\tilde{v}^{^{\mathrm{T}}}h \\
& -2\gamma_{2}\left( C\tilde{x} \right)^{\mathrm{T}} \left( R\otimes I_{q} \right) \tilde{v} -\gamma_{2} \tilde{v}^{\mathrm{T}} \left( R\otimes I_{q} \right) \tilde{v} \\
& + \tilde{v}^{\mathrm{T}} \left( R\otimes I_{q} \right) CA_{c}\left( \tilde{x}-\Psi\tilde{\rho} \right)\\
& + \delta_{c}\left(\tilde{x}-\Psi \tilde{\rho} \right)^{\mathrm{T}} \big( P_{c}A_{c}+ A_{c}^{\mathrm{T}}P_{c}\big)  \left(\tilde{x}-\Psi \tilde{\rho} \right)\\
& + 2 \left( C\tilde{x} \right)^{\mathrm{T}} \left( R\otimes I_{q} \right) CBK\breve{x} + \tilde{v}^{\mathrm{T}} \left( R\otimes I_{q} \right) CBK\breve{x}.
\end{aligned}
\end{equation*}
Then according to similar arguments in the proof of Theorem \ref{Theorem 1}, one obtains
\begin{align} \label{observer based V_1}
\dot{V}\!\leq\! & -\!\Big(\! 2m-\dfrac{M^{2}+\|C\|^{2}}{\delta} \Big) \lVert C\tilde{x} \rVert^{2} \notag\\
& -\!\Big(\!\gamma_{2} r_{\min}-\dfrac{5\delta}{4}- \dfrac{\|C\|^{2}}{4\delta} \Big) \|\tilde{v}\|^{2}   \notag\\
& -\!\Big(\! \gamma_{1} \lambda_{2}-\dfrac{\gamma_{2}^{2}}{\delta} \Big) \|s\|^{2} -\|\tilde{x}-\Psi \tilde{\rho}\|^{2} \notag\\
& +\! 2 \left( C\tilde{x} \right)^{\mathrm{T}} \!\left( R\otimes I_{q} \right) CBK\breve{x} \!+\! \tilde{v}^{\mathrm{T}}\! \left( R\otimes I_{q} \right) CBK\breve{x}. 
\end{align}

Note that the following inequalities are satisfied,
\begin{align} 
2 \left( C\tilde{x} \right)^{\mathrm{T}} &\left( R\otimes I_{q} \right) CBK\breve{x} \notag\\ 
& \leq \dfrac{\|C\|^{2}}{\delta} \|C\tilde{x}\|^{2} + \delta\|BK\|^{2} \|\breve{x}\|^{2}, \label{observer based inequality_1}\\
\tilde{v}^{\mathrm{T}} \big( R\otimes & I_{q} \big) CBK\breve{x} \notag\\ 
& \leq \dfrac{\|C\|^{2}}{4\delta} \|\tilde{v}\|^{2} + \delta\|BK\|^{2} \|\breve{x}\|^{2}. \label{observer based inequality_2}
\end{align}
Then substituting (\ref{observer based inequality_1}) and (\ref{observer based inequality_2}) into (\ref{observer based V_1}) leads to
\begin{align}  \label{observer based V_1.1}
\dot{V}\leq & -\Big( 2m-\dfrac{M^{2}+2\|C\|^{2}}{\delta} \Big) \lVert C\tilde{x} \rVert^{2} \notag\\
& -\Big(\gamma_{2} r_{\min}-\dfrac{5\delta}{4}- \dfrac{\|C\|^{2}}{2\delta} \Big) \|\tilde{v}\|^{2}   \notag\\
& -\Big( \gamma_{1} \lambda_{2}-\dfrac{\gamma_{2}^{2}}{\delta} \Big) \|s\|^{2} -\|\tilde{x}-\Psi \tilde{\rho}\|^{2} \notag\\
& + 2 \delta\|BK\|^{2} \|\breve{x}\|^{2}. 
\end{align}

Since $ A_{o} $ is Hurwitz, there exists a positive definite matrix $ P_{o} $ such that $ P_{o}A_{o}+ A_{o}^{\mathrm{T}}P_{o} \leq -I $. Consider another positive definite function $ V_{3}= \breve{x}^{\mathrm{T}}P_{o}\breve{x} $. The derivative of $ V_{3} $ along the trajectory of (\ref{observer based change system}) satisfies
\begin{equation*}
\begin{aligned}
\dot{V_{3}} = \breve{x}^{\mathrm{T}}(P_{o}A_{o}+ A_{o}^{\mathrm{T}}P_{o})\breve{x} \leq -\|\breve{x}\|^{2}.
\end{aligned}
\end{equation*}

Consider the Lyapunov function candidate $ \bar{V}(\tilde{x},\breve{x},\tilde{\rho},\tilde{v})= V_{1} + \delta_{c} V_{2} + \delta_{o} V_{3} $ with $ \delta_{o}= 2\delta \|BK\|^{2} +1 $, then the derivative of $ \bar{V} $ along the trajectory of (\ref{observer based change system}) satisfies
\begin{equation}
\begin{aligned} \notag
\dot{\bar{V}}\leq & -\Big( 2m-\dfrac{M^{2}+2\|C\|^{2}}{\delta} \Big) \lVert C\tilde{x} \rVert^{2} \\
& -\Big(\gamma_{2} r_{\min}-\dfrac{5\delta}{4}- \dfrac{\|C\|^{2}}{2\delta} \Big) \|\tilde{v}\|^{2}   \\
& -\Big( \gamma_{1} \lambda_{2}-\dfrac{\gamma_{2}^{2}}{\delta} \Big) \|s\|^{2}  -\|\tilde{x}-\Psi \tilde{\rho}\|^{2} -\|\breve{x}\|^{2} .  
\end{aligned}
\end{equation}

The rest of the proof follows arguments similar to those in the proof of Theorem \ref{Theorem 1}, and is thus omitted. \qed

\end{appendix}

\bibliographystyle{IEEEtran}  
\bibliography{IEEEabrv,mylib}

\begin{thebibliography}{10}
\providecommand{\url}[1]{#1}
\csname url@samestyle\endcsname
\providecommand{\newblock}{\relax}
\providecommand{\bibinfo}[2]{#2}
\providecommand{\BIBentrySTDinterwordspacing}{\spaceskip=0pt\relax}
\providecommand{\BIBentryALTinterwordstretchfactor}{4}
\providecommand{\BIBentryALTinterwordspacing}{\spaceskip=\fontdimen2\font plus
\BIBentryALTinterwordstretchfactor\fontdimen3\font minus
  \fontdimen4\font\relax}
\providecommand{\BIBforeignlanguage}[2]{{%
\expandafter\ifx\csname l@#1\endcsname\relax
\typeout{** WARNING: IEEEtran.bst: No hyphenation pattern has been}%
\typeout{** loaded for the language `#1'. Using the pattern for}%
\typeout{** the default language instead.}%
\else
\language=\csname l@#1\endcsname
\fi
#2}}
\providecommand{\BIBdecl}{\relax}
\BIBdecl

\bibitem{molzahn2017survey}
D.~K. Molzahn, F.~D{\"{o}}rfler, H.~Sandberg, S.~H. Low, S.~Chakrabarti,
  R.~Baldick, and J.~Lavaei, ``{A survey of distributed optimization and
  control algorithms for electric power systems},'' \emph{IEEE Transactions on
  Smart Grid}, vol.~8, no.~6, pp. 2941--2962, 2017.

\bibitem{ram2010distributed}
S.~S. Ram, V.~V. Veeravalli, and A.~Nedi{\'{c}}, ``{Distributed and recursive
  parameter estimation in parametrized linear state-space models},'' \emph{IEEE
  Transactions on Automatic Control}, vol.~55, no.~2, pp. 488--492, 2010.

\bibitem{tsitsiklis1984problems}
J.~N. Tsitsiklis, ``Problems in decentralized decision making and
  computation,'' Ph.D. dissertation, MIT, Cambridge, 1984.

\bibitem{nedic2009distributed}
A.~Nedi{\'{c}} and A.~Ozdaglar, ``{Distributed subgradient methods for
  multi-agent optimization},'' \emph{IEEE Transactions on Automatic Control},
  vol.~54, no.~1, pp. 48--61, 2009.

\bibitem{yang2019survey}
T.~Yang, X.~Yi, J.~Wu, Y.~Yuan, D.~Wu, Z.~Meng, Y.~Hong, H.~Wang, Z.~Lin, and
  K.~H. Johansson, ``{A survey of distributed optimization},'' \emph{Annual
  Reviews in Control}, vol.~47, pp. 278--305, 2019.

\bibitem{nedic2018distributed}
A.~Nedi{\'{c}} and J.~Liu, ``{Distributed optimization for control},''
  \emph{Annual Review of Control, Robotics, and Autonomous Systems}, vol.~1,
  no.~1, pp. 77--103, 2018.

\bibitem{liang2019distributed}
S.~Liang, L.~Y. Wang, and G.~Yin, ``{Distributed dual subgradient algorithms
  with iterate-averaging feedback for convex optimization with coupled
  constraints},'' \emph{IEEE Transactions on Cybernetics}, to be published,
  doi:
  \href{https://ieeexplore.ieee.org/abstract/document/8830450}{10.1109/TCYB.2019.2933003}.

\bibitem{nedic2014distributed}
A.~Nedi{\'c} and A.~Olshevsky, ``Distributed optimization over time-varying
  directed graphs,'' \emph{IEEE Transactions on Automatic Control}, vol.~60,
  no.~3, pp. 601--615, 2014.

\bibitem{xi2017dextra}
C.~Xi and U.~A. Khan, ``{DEXTRA: A fast algorithm for optimization over
  directed graphs},'' \emph{IEEE Transactions on Automatic Control}, vol.~62,
  no.~10, pp. 4980--4993, 2017.

\bibitem{xi2016distributed}
------, ``{Distributed subgradient projection algorithm over directed
  graphs},'' \emph{IEEE Transactions on Automatic Control}, vol.~62, no.~8, pp.
  3986--3992, 2016.

\bibitem{xi2018linear}
C.~Xi, V.~S. Mai, R.~Xin, E.~H. Abed, and U.~A. Khan, ``{Linear convergence in
  optimization over directed graphs with row-stochastic matrices},'' \emph{IEEE
  Transactions on Automatic Control}, vol.~63, no.~10, pp. 3558--3565, 2018.

\bibitem{wang2019distributed}
D.~Wang, J.~Yin, and W.~Wang, ``Distributed randomized gradient-free
  optimization protocol of multiagent systems over weight-unbalanced
  digraphs,'' \emph{IEEE Transactions on Cybernetics}, vol.~51, no.~1, pp.
  473--482, 2019.

\bibitem{wang2010control}
J.~Wang and N.~Elia, ``Control approach to distributed optimization,'' in
  \emph{Proc. 48th IEEE Annual Allerton Conference on Communication, Control,
  and Computing (Allerton)}, 2010, pp. 557--561.

\bibitem{gharesifard2013distributed}
B.~Gharesifard and J.~Cortes, ``{Distributed continuous-time convex
  optimization on weight-balanced digraphs},'' \emph{IEEE Transactions on
  Automatic Control}, vol.~59, no.~3, pp. 781--786, 2013.

\bibitem{kia2015distributed}
S.~S. Kia, J.~Cort{\'{e}}s, and S.~Mart{\'{i}}nez, ``{Distributed convex
  optimization via continuous-time coordination algorithms with discrete-time
  communication},'' \emph{Automatica}, vol.~55, pp. 254--264, 2015.

\bibitem{touri2016saddle}
B.~Touri and B.~Gharesifard, ``Saddle-point dynamics for distributed convex
  optimization on general directed graphs,'' in \emph{Proc. 55th IEEE
  Conference on Decision and Control (CDC)}, 2016, pp. 862--866.

\bibitem{li2017distributed}
Z.~Li, Z.~Ding, J.~Sun, and Z.~Li, ``{Distributed adaptive convex optimization
  on directed graphs via continuous-time algorithms},'' \emph{IEEE Transactions
  on Automatic Control}, vol.~63, no.~5, pp. 1434--1441, 2017.

\bibitem{zhu2018continuous}
Y.~Zhu, W.~Yu, G.~Wen, and W.~Ren, ``{Continuous-time coordination algorithm
  for distributed convex optimization over weight-unbalanced directed
  networks},'' \emph{IEEE Transactions on Circuits and Systems II: Express
  Briefs}, vol.~66, no.~7, pp. 1202--1206, 2018.

\bibitem{stegink2016Unifying}
T.~Stegink, C.~De~Persis, and A.~van~der Schaft, ``A unifying energy-based
  approach to stability of power grids with market dynamics,'' \emph{IEEE
  Transactions on Automatic Control}, vol.~62, no.~6, pp. 2612--2622, 2016.

\bibitem{song2017Relative}
W.~Song, Y.~Tang, Y.~Hong, and X.~Hu, ``{Relative attitude formation control of
  multi-agent systems},'' \emph{International Journal of Robust and Nonlinear
  Control}, vol.~27, no.~18, pp. 4457--4477, 2017.

\bibitem{zhang2011extremum}
C.~Zhang and R.~Ord{\'o}{\~n}ez, \emph{Extremum-seeking control and
  applications: a numerical optimization-based approach}.\hskip 1em plus 0.5em
  minus 0.4em\relax Springer Science \& Business Media, 2011.

\bibitem{tran2018distributed}
N.~T. Tran, Y.~W. Wang, and W.~Yang, ``{Distributed optimization problem for
  double-integrator systems with the presence of the exogenous disturbance},''
  \emph{Neurocomputing}, vol. 272, pp. 386--395, 2018.

\bibitem{xie2019global}
Y.~Xie and Z.~Lin, ``{Global optimal consensus for higher-order multi-agent
  systems with bounded controls},'' \emph{Automatica}, vol.~99, pp. 301--307,
  2019.

\bibitem{tang2018optimal}
Y.~Tang, Z.~Deng, and Y.~Hong, ``{Optimal output consensus of high-order
  multi-agent systems with embedded technique},'' \emph{IEEE Transactions on
  Cybernetics}, vol.~49, no.~5, pp. 1768--1779, 2018.

\bibitem{wang2015distributed}
X.~Wang, Y.~Hong, and H.~Ji, ``{Distributed optimization for a class of
  nonlinear multi-agent systems with disturbance rejection},'' \emph{IEEE
  Transactions on Cybernetics}, vol.~46, no.~7, pp. 1655--1666, 2015.

\bibitem{tang2018distributed}
Y.~Tang, ``{Distributed optimization for a class of high-order nonlinear
  multi-agent systems with unknown dynamics},'' \emph{International Journal of
  Robust and Nonlinear Control}, vol.~28, no.~17, pp. 5545--5556, 2018.

\bibitem{zhao2017distributed}
Y.~Zhao, Y.~Liu, G.~Wen, and G.~Chen, ``{Distributed optimization for linear
  multi-agent systems: edge-and node-based adaptive designs},'' \emph{IEEE
  Transactions on Automatic Control}, vol.~62, no.~7, pp. 3602--3609, 2017.

\bibitem{li2019distributed}
Z.~Li, Z.~Wu, Z.~Li, and Z.~Ding, ``{Distributed optimal coordination for
  heterogeneous linear multi-agent systems with event-triggered mechanisms},''
  \emph{IEEE Transactions on Automatic Control}, vol.~65, no.~4, pp.
  1763--1770, 2019.

\bibitem{zhang2020exponential}
J.~Zhang, L.~Liu, and H.~Ji, ``Exponential convergence of distributed optimal
  coordination for linear multi-agent systems over general digraphs,'' in
  \emph{Proc. 39th Chinese Control Conference (CCC)}, 2020, pp. 5047--5051.

\bibitem{huang2004nonlinear}
J.~Huang, \emph{Nonlinear output regulation: theory and applications}.\hskip
  1em plus 0.5em minus 0.4em\relax SIAM, 2004.

\bibitem{bullo2019lectures}
F.~Bullo, \emph{Lectures on network systems}.\hskip 1em plus 0.5em minus
  0.4em\relax Kindle Direct Publishing, 2019.

\bibitem{olfati-Saber2004consensus}
R.~Olfati-Saber and R.~M. Murray, ``{Consensus problems in networks of agents
  with switching topology and time-delays},'' \emph{IEEE Transactions on
  Automatic Control}, vol.~49, no.~9, pp. 1520--1533, 2004.

\bibitem{ren2005consensus}
W.~Ren and R.~W. Beard, ``{Consensus seeking in multi-agent systems under
  dynamically changing interaction topologies},'' \emph{IEEE Transactions on
  Automatic Control}, vol.~50, no.~5, pp. 655--661, 2005.

\bibitem{bertsekas2009convex}
D.~P. Bertsekas, \emph{Convex optimization theory}.\hskip 1em plus 0.5em minus
  0.4em\relax Athena Scientific Belmont, 2009.

\bibitem{bernstein2009matrix}
D.~S. Bernstein, \emph{Matrix mathematics: theory, facts, and formulas},
  2nd~ed.\hskip 1em plus 0.5em minus 0.4em\relax Princeton University Press,
  2009.

\bibitem{khalil2002nonlinear}
H.~K. Khalil, \emph{Nonlinear systems}, 3rd~ed.\hskip 1em plus 0.5em minus
  0.4em\relax Prentice Hall, 2002.

\bibitem{wieland2011internal}
P.~Wieland, R.~Sepulchre, and F.~Allg{\"o}wer, ``An internal model principle is
  necessary and sufficient for linear output synchronization,''
  \emph{Automatica}, vol.~47, no.~5, pp. 1068--1074, 2011.

\bibitem{kim2010output}
H.~Kim, H.~Shim, and J.~H. Seo, ``Output consensus of heterogeneous uncertain
  linear multi-agent systems,'' \emph{IEEE Transactions on Automatic Control},
  vol.~56, no.~1, pp. 200--206, 2010.

\bibitem{chen1999linear}
C.-T. Chen, \emph{Linear system theory and design}, 3rd~ed.\hskip 1em plus
  0.5em minus 0.4em\relax Oxford University Press, Inc., 1999.

\end{thebibliography}

\end{document}